\documentclass{article}

\usepackage{amsmath}
\usepackage{amssymb}
\usepackage{commath}

\usepackage{float}
\usepackage{multirow}

\usepackage{xcolor}

\usepackage{titlesec}
\usepackage[section]{placeins}

\usepackage[utf8]{inputenc}
\usepackage[english]{babel}

\usepackage{etoolbox}

\usepackage{graphicx}
    \graphicspath{{img/}{../img/}}
\usepackage[labelformat=simple]{subcaption}

\usepackage{indentfirst}
\usepackage{enumitem}

\usepackage{amsmath}
\usepackage{amssymb}
\usepackage{commath}

\usepackage{float}
\usepackage{hyperref}
\usepackage[all]{hypcap}
\usepackage{multirow}

\usepackage{xcolor}

\usepackage{titlesec}

\newcommand{\me}{\mathrm{e}}  
\newcommand{\iu}{{i\mkern1mu}} 
\newcommand{\rpm}{\sbox0{$1$}\sbox2{$\scriptstyle\pm$} 
  \raise\dimexpr(\ht0-\ht2)/2\relax\box2 }
  
\makeatletter 
\renewenvironment{cases}[1][l]{\matrix@check\cases\env@cases{#1}}{\endarray\right.}
\def\env@cases#1{%
  \let\@ifnextchar\new@ifnextchar
  \left\lbrace\def\arraystretch{1.2}%
  \array{@{}#1@{\quad}l@{}}}
\makeatother

\usepackage{etoolbox}
\pretocmd{\eqref}{Eq.~}{}{}

\newcommand{\Secref}[1]{Section~\ref{#1}} 
\newcommand{\tableref}[1]{Table~\ref{#1}}
\newcommand{\listref}[1]{Listing~\ref{#1}}
\newcommand{\R}{\mathbb{R}}

\usepackage{listings}
\usepackage{authblk}

\title{Interactive $G^1$ and $G^2$ Hermite Interpolation\\ Using Extended Log-aesthetic Curves}

\author[1,2]{Ferenc Nagy}
\author[3]{Norimasa Yoshida}
\author[1,4]{Mikl\'{o}s Hoffmann}

\affil[1]{University of Debrecen, Faculty of Informatics}
\affil[2]{University of Debrecen, Doctoral School of Informatics}
\affil[3]{Nihon University, College of Industrial Technology}
\affil[4]{Eszterh\'azy K\'aroly University, Institute of Mathematics and Computer Science}

\date{March 2021}

\begin{document}

\maketitle

\abstract{In the field of aesthetic design, log-aesthetic curves have a significant role to meet the high industrial requirements. In this paper, we propose a new interactive $G^1$ Hermite interpolation method based on the algorithm of Yoshida et al. \cite{yoshida2006interactive} with a minor boundary condition. In this novel approach, we compute an extended log-aesthetic curve segment that may include inflection point (S-shaped curve) or cusp. The curve segment is defined by its endpoints, a tangent vector at the first point, and a tangent direction at the second point. The algorithm also determines the shape parameter of the log-aesthetic curve based on the length of the first tangent that provides control over the curvature of the first point and makes the method capable of joining log-aesthetic curve segments with $G^2$ continuity.
}

\section{Introduction}

Aesthetic curves are primarily used in computer-aided design to meet the high aesthetic requirements of the industry. Levien et al. stated \cite{levien2009interpolating} that the log-aesthetic curve is the most promising curve for aesthetic design and a large number of research papers are published since their introduction.

The log-aesthetic curve is originated from Harada et al. \cite{harada1997study, harada1999aesthetic}. They analyzed the characteristics of aesthetic curves, and insisted that natural aesthetic curves have such a property that their logarithmic distribution diagram of curvature (LDDC) can be approximated by straight lines meanwhile there is a strong correlation between the slopes of the lines and the impressions of the curves.
Based on their work Miura et al. \cite{miura2005derivation,miura2006general} have defined the logarithmic curvature graph (LCG), an analytical version of the LDDC as follows: when a curve (given with arc length $s$ and radius of curvature $\rho$) is subdivided into infinitesimal segments such that $\Delta\rho/\rho$ is constant, the LCG represents the relationship between $\rho$ and $\Delta s$ in a double logarithmic graph.
They have also formulated the curve whose LCG is strictly expressed by a straight line (with slope $\alpha$):
\begin{equation} \label{eq:fundamental}
	\log\rho\frac{\dif s}{\dif \rho} = \alpha \log \rho + c,
\end{equation}
where $c$ is a constant.
\eqref{eq:fundamental} is the fundamental equation of log-aesthetic curves.

Yoshida et al. \cite{yoshida2006interactive} analyzed the properties of the log-aesthetic curve and derived a general formula from the relationship between the arc length and the radius of curvature of the curve. The authors are also presented an interactive algorithm to draw a log-aesthetic curve segment, which is an essential approach in log-aesthetic design. Several algorithms are developed based on their technique. 

However, their method can only generate a log-aesthetic curve segment with monotonic curvature variation and can not create a curve with curvature-extremal point or inflection point. Miura et al. \cite{miura2008input} presented a novel technique to input a log-aesthetic curve segment with an inflection point and a method to generate a log-aesthetic curve from a sequence of 2D points.
However, their technique restricts users from specifying tangent directions at the endpoints. Hence, the proposed method cannot solve $G^2$ Hermite interpolation problem.
Therefore, in \cite{Miura2013} Miura et al. proposed a new method to generate an S-shaped log-aesthetic curve and a novel method to solve the $G^2$ Hermite interpolation problem with log-aesthetic curves in the form of log-aesthetic triplets.
Besides, Meek et al. \cite{meek2012} used planar log-aesthetic spirals that include a point of zero curvature and proved that for any member of the family a unique segment of that spiral can be found that matches the given two-point $G^1$ Hermite data. They investigated the cases when the shape parameter $\alpha < 1$ holds. Furthermore, a so-called generalized log-aesthetic curve has been developed \cite{gobithaasan2011aesthetic} which extends the log-aesthetic curve by expressing LCG in a linear form. The generalized log-aesthetic curve provides the possibility of control of its curvature with an extra shape parameter.

In the algorithms dealing with log-aesthetic curve, there is usually some restrictions on the curve, e.g. there is a problem that computation may not be possible for the given boundary condition depending on the shape parameter \cite{yoshida2012drawable,gobithaasan2013drawable}. A potential solution is the discretization of the log-aesthetic curve, proposed by Yagi et al. \cite{yagi2020g1}.

In this paper, we propose a new, interactive $G^1$ Hermite interpolation method based on the fundamental algorithm of Yoshida et al. \cite{yoshida2006interactive}. In this approach, the log-aesthetic curve may include inflection point, i.e. can form S-shaped curve, or it can include curvature-extremal point, on the user's desire. We call the new curve extended log-aesthetic curve. The extension provides a solution to design aesthetic curves through geometric data with a minor boundary condition.

The new method generates an extended log-aesthetic curve segment defined by two points, a tangent vector at the first point, and a tangent line at the last point. This method is much closer to the classical Hermite interpolation method (based on two points and two tangent vectors) than the previous log-aesthetic methods. The $\alpha$ parameter is determined based on the length of the first tangent vector. 
In the new algorithm the user also possesses control over the curvature of the first point that makes the method capable of joining log-aesthetic curve segments with $G^2$ continuity, because the tangent vector of the last point can be calculated from the derivative of the tangential angle parametrized curve, that equals the radius of curvature $\rho(\theta)$
(further discussed in \Secref{sec:previous_results}).
Therefore, if the lengths of the tangent vectors match (besides the direction), the log-aesthetic curve segments are also share a common curvature at their joint.
We must remark that in our algorithm we exclude the very specific case, when the given tangent vector of the first point is parallel to the tangent line of the last point.  

Since the new approach relies on the work of Yoshida et al. \cite{yoshida2006interactive}, their derived general formulae of the log-aesthetic curve and their previous algorithm will be briefly described in \Secref{sec:previous_results}.
The extended log-aesthetic curves that can include inflection point and cusp will be presented in \Secref{sec:extended_log-aesthetic_curves}.
The new method will be described in detail in \Secref{sec:lambda_bisection} and in \Secref{sec:alpha_bisection}.

\section{Previous Results} \label{sec:previous_results}

In this section, we are going to introduce the notations and summarise the results of \cite{yoshida2006interactive}. We also briefly describe the primary interactive modeling algorithm.

The general formula is derived using a reference point (at the origin) and placing the following constraints on it: the radius of curvature is 1 and the tangent vector is directed in the positive direction along the x coordinate axis. The standard form is obtained by transforming the log-aesthetic curve such that the above constraints are satisfied, where $\Lambda$ is responsible for the transformations.
Therefore, the point $P(\theta)$ of a log-aesthetic curve whose tangential angle is $\theta$, is expressed on the complex plane as the following (\cite{yoshida2006interactive}):
\begin{equation} \label{eq:log-aesthetic_point_theta}
    P(\theta)= 
\begin{cases}[c]
    \int_{0}^{\theta}  \me^{(1+\iu)\Lambda\psi}\  d\psi
        & \text{if } \alpha = 1\\
    \int_{0}^{\theta}  \big((\alpha-1)\Lambda\psi+1\big)^{\frac{1}{\alpha-1}} \me^{\iu\psi}\  d\psi
        & \text{otherwise},
\end{cases}
\end{equation}
where $\alpha \in \R$ and $\Lambda \in \R^{+}$ are parameters. $\alpha$ is the slope of the LCG.
The point of $\theta=0$ is the reference point at the origin and its tangent vector is $\big[1$ $0\big]^T$. The $\Lambda$ corresponds to the constant part of \eqref{eq:fundamental}.
When $\alpha \neq 1$, all the aesthetic curves are congruent under similarity transformations without depending on the value of $\Lambda(\neq 0$).

The radius of curvature of the log-aesthetic curve is (\cite{yoshida2006interactive}):
\begin{equation} \label{eq:radius_curvature_theta}
    \rho(\theta)= 
\begin{cases}[c]
    \me^{\Lambda\theta}
        & \text{if } \alpha = 1\\
    \big((\alpha-1)\Lambda\theta+1\big)^{\frac{1}{\alpha-1}} 
        & \text{otherwise}.
\end{cases}
\end{equation}
When $\theta=0$, $\rho=1$. The function increases monotonically when $\Lambda \neq 0$. In case of $\Lambda = 0$, $\rho$ is constant $1$ and the log-aesthetic curve is a circle. Similarly, when $\alpha = \rpm \infty$ because taking the limit of \eqref{eq:radius_curvature_theta} when $\alpha$ approaches $\rpm \infty$ we get $\rho=1$.

The tangential angle $\theta$ and the arc length $s$ are related by (\cite{yoshida2006interactive}):
\begin{equation} \label{eq:theta_arc}
    \theta(s) = 
\begin{cases}[c]
    \frac{1-\me^{-\Lambda s}}{\Lambda}  
        & \text{if } \alpha = 0\\
    \frac{\log(\Lambda s + 1 )}{\Lambda}
        & \text{if } \alpha = 1\\
    \frac{(\Lambda\alpha s +1)^{(1-\frac{1}{\alpha})}-1} {\Lambda(\alpha-1)}  
        & \text{otherwise}.
\end{cases}
\end{equation}
Therefore, a point on the aesthetic curve $C(s)$ whose arc length is $s$, can be defined on the complex plane as (\cite{yoshida2006interactive}):
\begin{equation} \label{eq:log-aesthetic_point_arc}
    C(s)= 
\begin{cases}[c]
    \int_{0}^{s}\exp({\iu\frac{1-\me^{-\Lambda u}}{\Lambda}})\ d u
        & \text{if } \alpha = 0\\
    \int_{0}^{s}\exp({\iu\frac{\log(\Lambda u + 1 )}{\Lambda}})\ d u
        & \text{if } \alpha = 1\\
    \int_{0}^{s}\exp({\iu\frac{(\Lambda\alpha u +1)^{(1-\frac{1}{\alpha})}-1} {\Lambda(\alpha-1)}})\  d u
        & \text{otherwise}.
\end{cases}
\end{equation}

\eqref{eq:log-aesthetic_point_arc} and \eqref{eq:log-aesthetic_point_theta} represent the same curve.
Using \eqref{eq:theta_arc}, the radius of curvature can also be expressed from the arc length $s$ \cite{yoshida2006interactive}:
\begin{equation} \label{eq:radius_curvature_arc}
    \rho(s)= 
\begin{cases}[c]
    \me^{\Lambda s}
        & \text{if } \alpha = 0\\
    \big(\Lambda\alpha s + 1 \big)^{\frac{1}{\alpha}} 
        & \text{otherwise.}
\end{cases}
\end{equation}

Since $\rho$ can change from $-\infty$ to $+\infty$, the tangential angle $\theta$ and arc length $s$ may have upper or lower bound depending on the value of $\alpha$ (because of the negative bases of the fractional exponents):
\begin{table}[H]
\centering
\begin{tabular}{|c c c c c c c c |}
\hline
    & \multicolumn{3}{c}{\textbf{Tangential angle ($\boldsymbol{\theta}$)}} &
    & \multicolumn{3}{c|}{\textbf{Arc length ($\boldsymbol{s}$)}} \\
\hline
    & $\alpha<1$ & $\alpha=1$ & $\alpha>1$ & &
      $\alpha<0$ & $\alpha=0$ & $\alpha>0$ \\
    Upper bound:  &
    $\frac{1}{\Lambda(1-\alpha)}$ & - & - & &
    $-\frac{1}{\Lambda\alpha}$ & - & - \\
    Lower bound:  &
    - & - & $\frac{1}{\Lambda(1-\alpha)}$ & &
    - & - & $-\frac{1}{\Lambda\alpha}$ \\
\hline
\end{tabular}
\caption[]{Upper and lower bound of $\theta$ and $s$.}
\label{table:arc_theta_bounds}
\end{table}

It is worth to mention that the derived formulas can be numerically unstable for $\alpha \approx 0$ or $\alpha \approx 1$ but $\alpha \neq 0$ and $\alpha \neq 1$. However, these cases can be avoided simply by ignoring the neighbouring $\alpha$ values $\left] 0, 0\rpm\epsilon\right]$ and $\left] 1, 1\rpm\epsilon\right]$ , which causes only insignificant deviation from the desired geometric data during the algorithm. \cite{meek2012}

Besides the equations above and the properties of the log-aesthetic curves, Yoshida et al. also presented an interactive algorithm \cite{yoshida2006interactive} to generate a log-aesthetic curve segment by specifying three so-called control points (similarly in case of a quadratic B\'ezier curve) and $\alpha$.
The idea is to search for a curve segment that fits a similar triangle defined by the control points, using a bisection method on $\Lambda$.  
The curve is drawn from the first control point $A$ to the last $C$, while point $B$ specifies the change of the tangential angle $\theta_\Delta$ between the endpoints. If $\left| A B \right| \leq \left| B C \right|$ does not hold, the coordinates of $A$ and $C$ need to be swapped.

\begin{figure}[!ht]
  \subfloat[$\alpha \leq 1$]{
	\begin{minipage}[t]{0.95\textwidth}
	\label{fig:previous_algorithm_conf_1}
	   \centering
	   \includegraphics[width=0.95\textwidth]{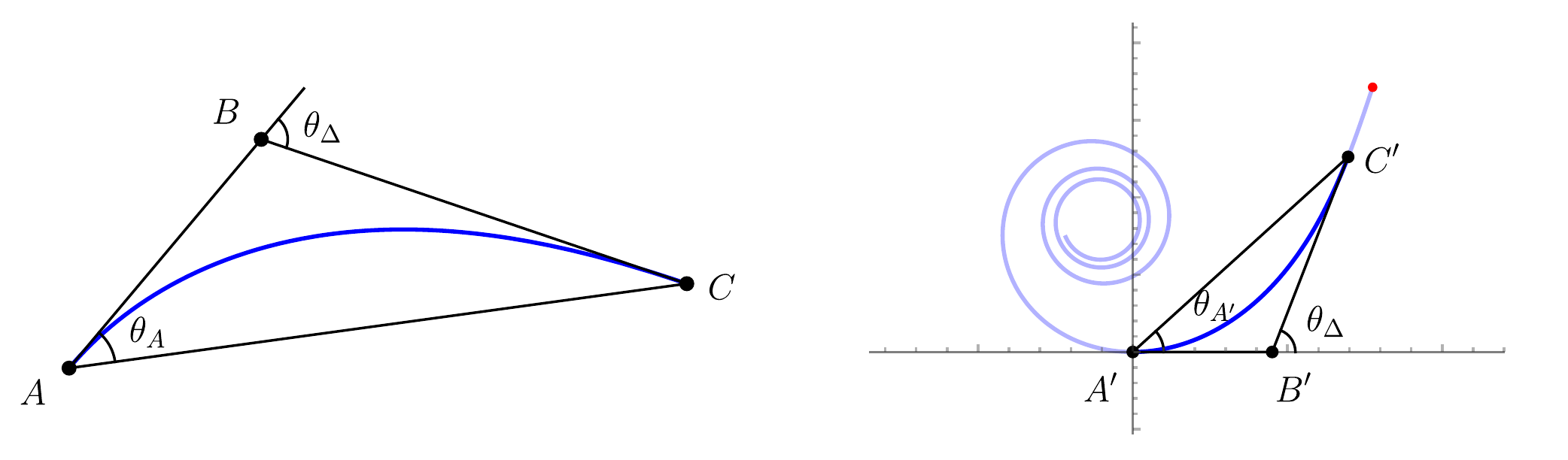}
	\end{minipage}}
	\\
  \subfloat[$\alpha > 1$]{
	\begin{minipage}[t]{0.95\textwidth}
	\label{fig:previous_algorithm_conf_2}
	   \centering
	   \includegraphics[width=0.95\textwidth]{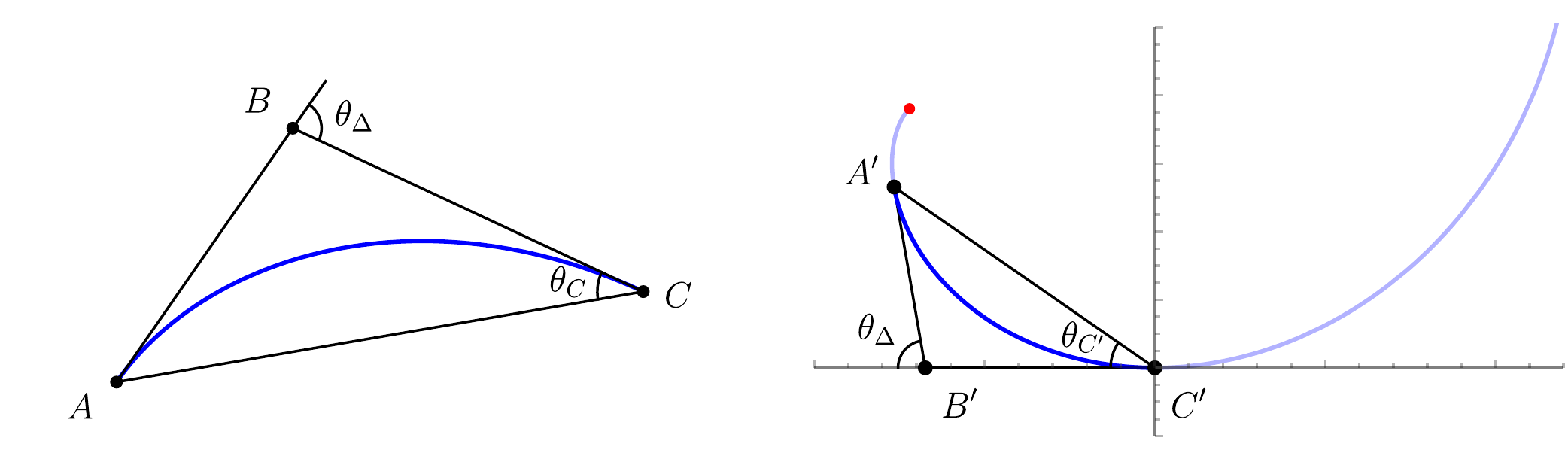}
	\end{minipage}}
\caption{Configurations of the previous algorithm, where the red point is the bound of the log-aesthetic curve (based of Fig. 7 in \cite{yoshida2006interactive})}
\label{fig:previous_algorithm}
\end{figure}

The algorithm of \cite{yoshida2006interactive} works as follows, referring to \figref{fig:previous_algorithm}. The point $A'$ that corresponds to $A$ is defined as the point of the log-aesthetic curve whose tangential angle is 0. The point $C'$ corresponds to $C$, where the tangential angle is $\theta_\Delta$. The point $B'$ is the intersection of the tangent lines of the endpoints. It is placed on the x-axis because the tangent vector of $A'$ is \mbox{$\big[1$ $0\big]^T$}. If the triangle $A' B' C'$ is similar to the triangle $ABC$, the curve segment can be drawn by transforming the log-aesthetic points between the triangle $A'B'C'$ to the triangle $ABC$ (See \figref{fig:previous_algorithm}). The algorithm uses a bisection method on $\Lambda$ to find a similar triangle. The similarity is tested by comparing the corresponding angles of the triangles: $\theta_A$ with $\theta_{A'}$ in case of $\alpha \leq 1$ or $\theta_C$ with $\theta_{C'}$ otherwise. If the angles are equal, the triangles are similar because the exterior angle at point $B$ and $B'$ is $\theta_\Delta$.

When $\alpha<1$, the integration range of $[0,\theta_\Delta]$ is used. In this case 
$\theta$ has upper bound so $\theta_\Delta < 1/(\Lambda(1-\alpha))$ must hold. Therefore $0 < \Lambda < 1/(\theta_\Delta(1-\alpha))$.
When $\alpha>1$, $A'$ is defined as the point where the tangential angle is $-\theta_\Delta$ and $C'$ is the point of $\theta=0$. In this case, the range of the integration is $[-\theta_\Delta,0]$ and since $\theta$ has lower bound $-\theta_\Delta > 1/(\Lambda(1-\alpha))$ must hold. Therefore $\Lambda$ is between 0 and $1/(\theta_\Delta(\alpha-1))$.
When $\alpha=1$, there is no upper or lower bound for $\theta$, the bisection method is extended so that $\Lambda(>0)$ can be arbitrarily large.
For the pseudo-code of the bisection method, see Appendix A in \cite{yoshida2006interactive}.

\begin{figure}[!ht]
  \centering
  \subfloat[$\Lambda = 0$]{
	\begin{minipage}[t][1\width]{0.23\textwidth}
	   \centering
	   \includegraphics[width=\textwidth]{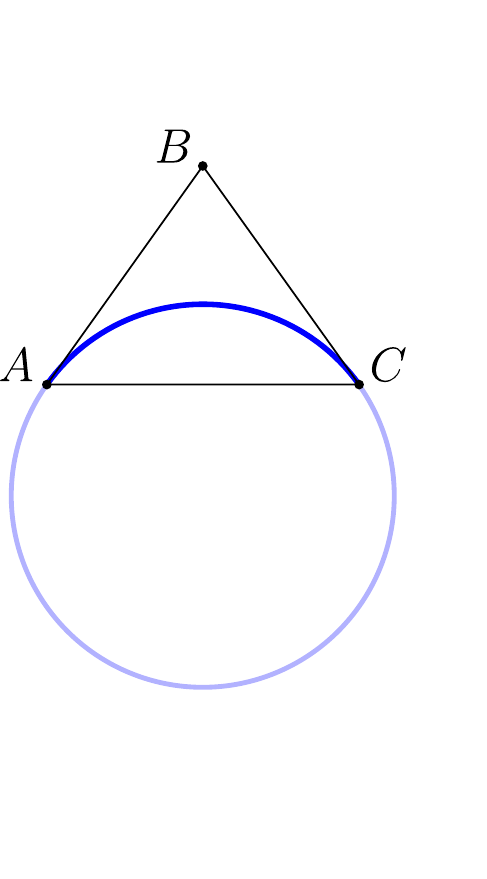}
	\end{minipage}}
 \hfill 	
  \subfloat[$\Lambda \approx 0.178$]{
	\begin{minipage}[t][1\width]{0.23\textwidth}
	   \centering
	   \includegraphics[width=\textwidth]{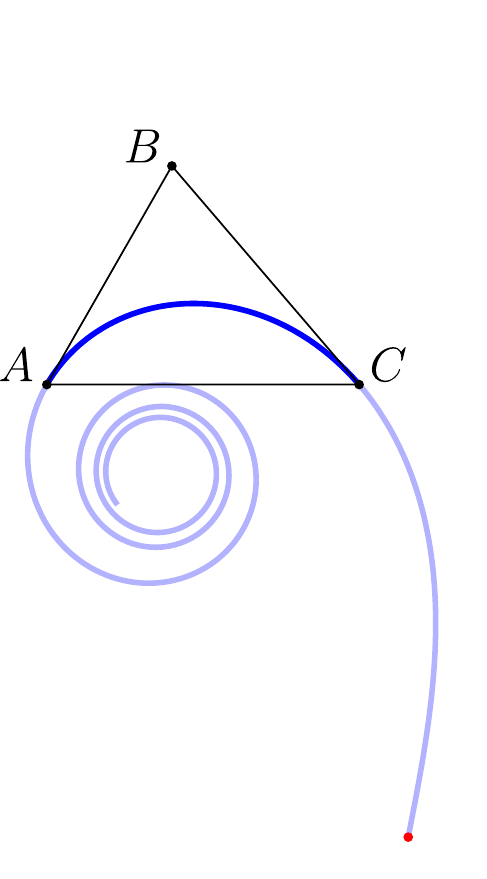}
	\end{minipage}}
  \hfill 	
  \subfloat[$\Lambda \approx 0.234$]{
	\begin{minipage}[t][1\width]{0.23\textwidth}
	   \centering
	   \includegraphics[width=\textwidth]{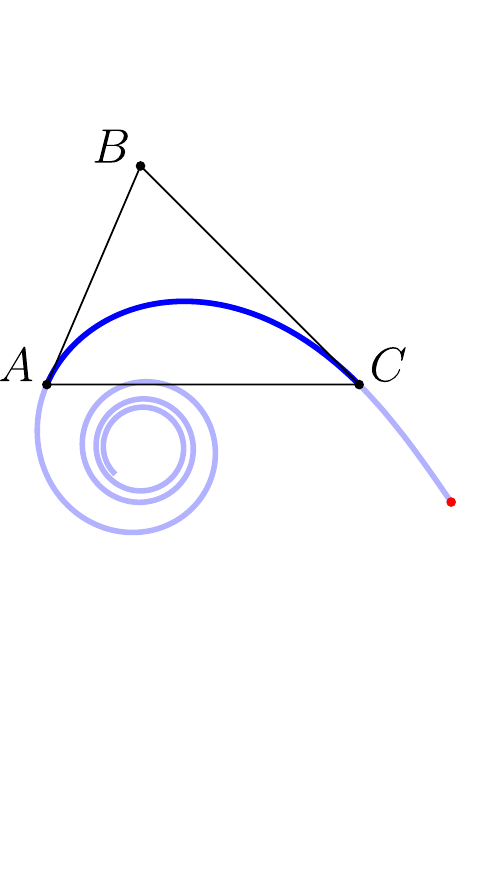}
	\end{minipage}}
  \hfill 	
  \subfloat[$\Lambda \approx 0.245$]{
	\begin{minipage}[t][1\width]{0.23\textwidth}
	   \centering
	   \includegraphics[width=\textwidth]{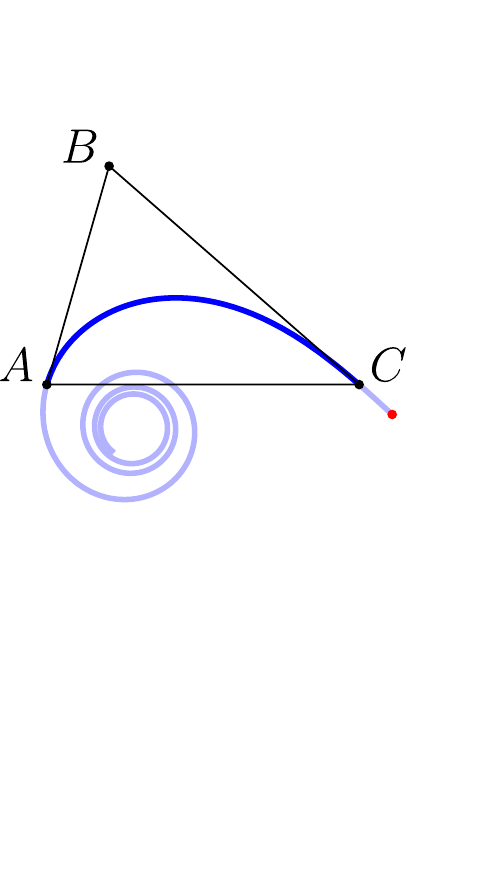}
	\end{minipage}}
\caption{Effect of alteration of the position of point $B$ ($\alpha=-1$). In case of further change, the bound happens to be inside the triangle and the algorithm has no solution.}
\label{fig:previous_algorithm_error}
\end{figure}

The main drawback of their algorithm is that the position of the control points and the value of $\alpha$ highly restrict the region where the curve can be drawn. Using the configuration above the largest log-aesthetic curve segment (based on the arc length $s$ between $A'$ and $C'$) can be drawn from the reference point up to the bound or vice versa. If the control points form an isosceles triangle such that the sides $AB$ and $BC$ are equal, the drawn log-aesthetic curve segment is a circular arc with $\Lambda=0$. If we move point $B$ parallel to the line $AC$ toward an endpoint, this endpoint will be the first point of the curve segment that corresponds to the reference point because of the coordinate swapping.
By this move, the determined $\Lambda$ increases and the other endpoint approaches the bound of the log-aesthetic curve (see \figref{fig:previous_algorithm_error}). The curve segment can be determined until the endpoint is within the bound.
In \cite{yoshida2006interactive}, the area where the curve segment can be determined is called the drawable region. If $\alpha$ is less than 0 or greater than 1, this drawable region gets drastically smaller (See e.g. Figure 9 in \cite{yoshida2006interactive}). For further survey on boundaries, see \cite{gobithaasan2013drawable} and \cite{yoshida2012drawable}.
In the followings, to extend the capability of their algorithm, we introduce the new, extended log-aesthetic curve.

\section{New Algorithm}

As it can be seen from the algorithm of \cite{yoshida2006interactive} the capability to use log-aesthetic curves is limited, due to the bounds. Our idea is to extend the curve beyond the bound by concatenating log-aesthetic segments to make the drawable regions larger.
In this section, we are going to propose a way for the extension and describe the new $\Lambda$ bisection method that works on the extended log-aesthetic curves.

The ultimate goal is to find an appropriate value for $\alpha$ to determine a log-aesthetic segment with the given tangent direction and length at the first point, in addition to the $G^1$ Hermite conditions. The last subsection will describe an alternative bisection method, applied to $\alpha$ to search for its appropriate value making the algorithm capable of joining log-aesthetic curves even with $G^2$ continuity.

\subsection{Extended Log-aesthetic Curves} \label{sec:extended_log-aesthetic_curves}

The extension of the log-aesthetic curve is defined by mirroring the curve at the original bounds. Regarding the value of $\alpha$, four cases can be distinguished: when $\alpha > 1$, $\alpha = 1$, $0 \leq \alpha < 1$, and $\alpha < 0$. By defining the reflection we take into consideration of the integration range of the previous algorithm since our aim is to extend its capability.

\begin{figure}[!ht]
  \centering
  \subfloat[$\alpha=1.5$]{
	\begin{minipage}[c][1\width]{0.45\textwidth}
	\label{fig:log-aesthetic_alpha_1_5}
	   \centering
	   \includegraphics[width=\textwidth]{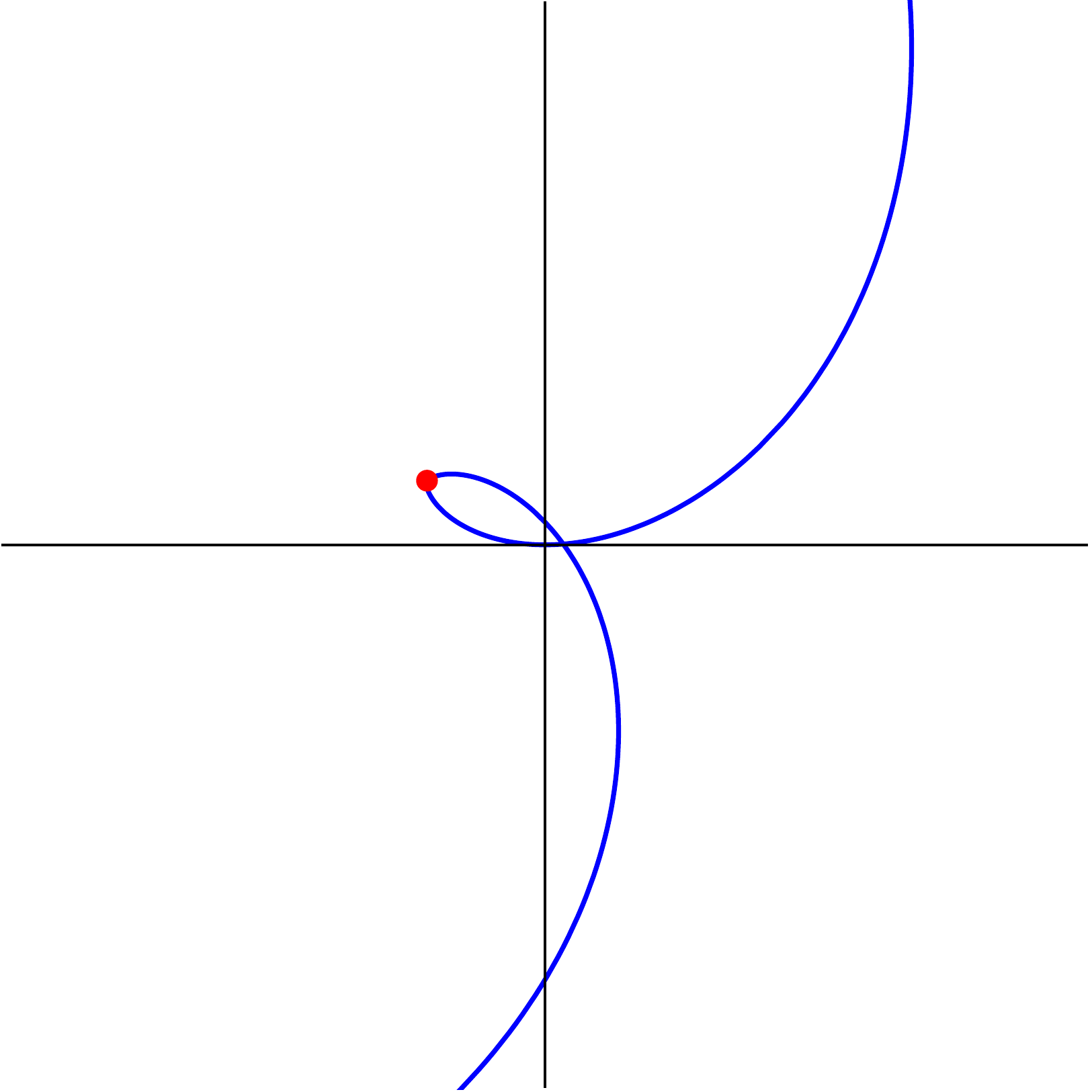}
	\end{minipage}}
 \hfill 	
  \subfloat[$\alpha=2$]{
	\begin{minipage}[c][1\width]{0.45\textwidth}
	\label{fig:log-aesthetic_alpha_2}
	   \centering
	   \includegraphics[width=\textwidth]{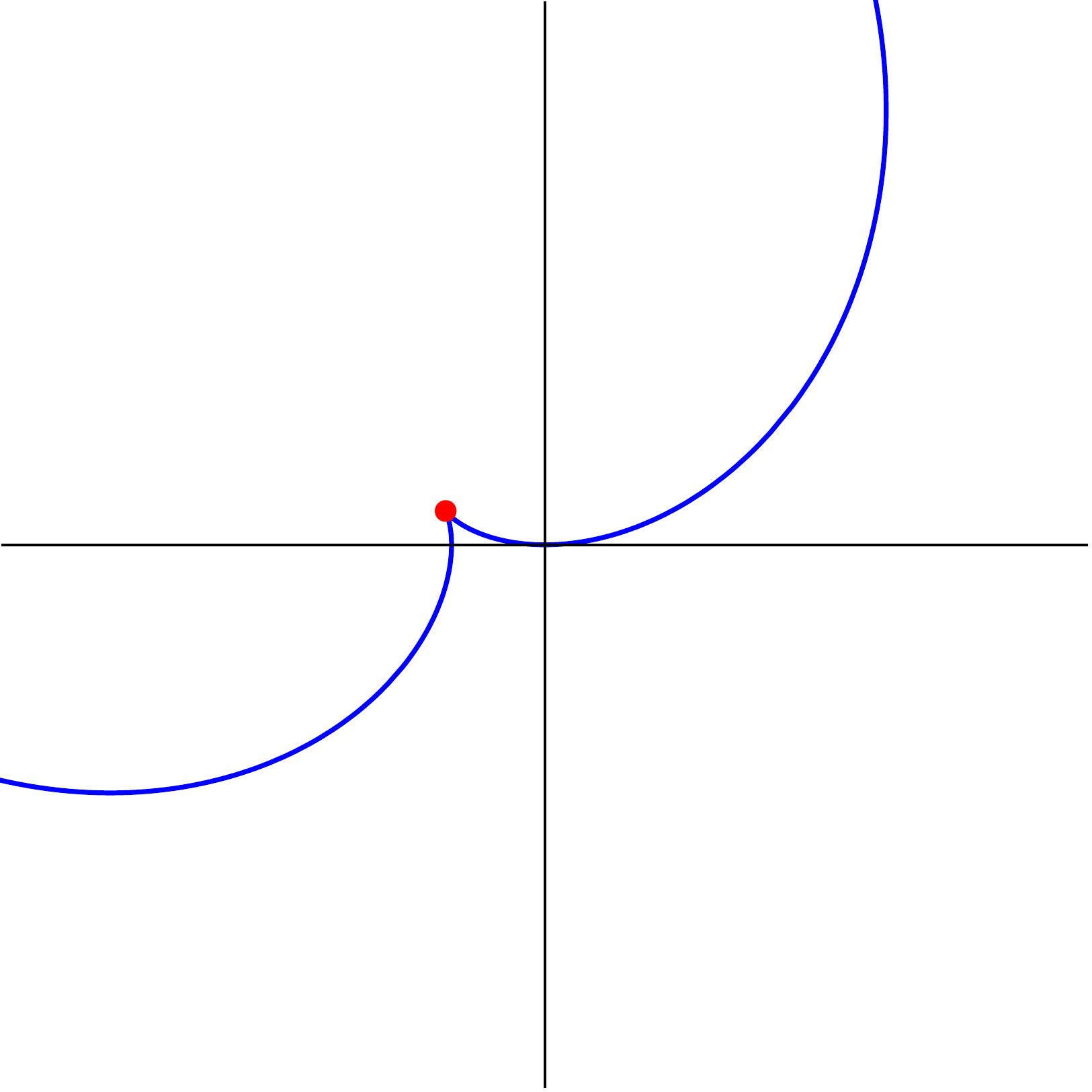}
	\end{minipage}}
\caption{Log-aesthetic curves, without applying the bound of $\theta$. We prefer the case of Figure (b).}
\end{figure}

In case of $\alpha>1$, the integration range of the above algorithm is $[0,-\theta_\Delta]$ but the curve has a lower bound at $\theta=\frac{1}{\Lambda(1-\alpha)}$. Therefore, we need to extend the curve beyond the bound to increase the drawable region. Since $\rho=0$ at the point of bound, the directional derivative is $0$ (see \eqref{eq:log-aesthetic_point_theta}). It is a singular point of the curve. We apply the bound of \tableref{table:arc_theta_bounds} on $\theta$ due to the possible negative bases of the fractional exponent. However, there are some $\alpha$ cases when the curve can be calculated on the entire domain of $\theta$.
For example when $\alpha=1.5$ (see \figref{fig:log-aesthetic_alpha_1_5})
the \eqref{eq:log-aesthetic_point_theta} is
\begin{equation} \label{eq:log-aesthetic_point_theta_alpha_1_5}
    P(\theta)= 
    \int_{0}^{\theta}  \left(\frac{\psi \Lambda }{2}+1\right)^2 \me^{\iu\psi}\  d\psi,
\end{equation}
or in case $\alpha=2$ (when the log-aesthetic curve is a circle involute, see \figref{fig:log-aesthetic_alpha_2}) the \eqref{eq:log-aesthetic_point_theta} is
\begin{equation} \label{eq:log-aesthetic_point_theta_alpha_2}
    P(\theta)= 
    \int_{0}^{\theta} \psi \Lambda +1 \me^{\iu\psi}\  d\psi.
\end{equation}

\begin{figure}[!ht]
  \centering
  \subfloat[$\alpha=1.5$]{
	\begin{minipage}[c][1\width]{0.45\textwidth}
	\label{fig:rho_alpha_1_5}
	   \centering
	   \includegraphics[width=\textwidth]{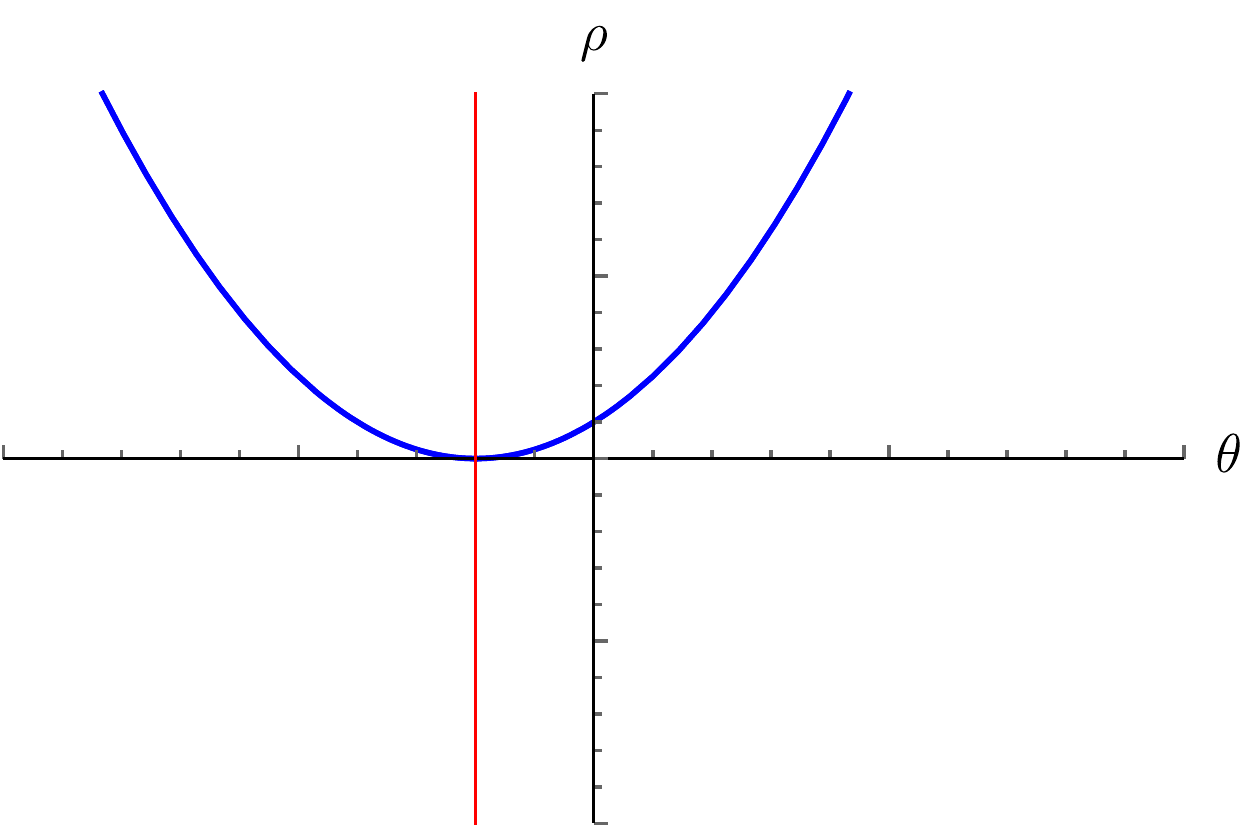}
	\end{minipage}}
 \hfill 	
  \subfloat[$\alpha=2$]{
	\begin{minipage}[c][1\width]{0.45\textwidth}
	\label{fig:rho_alpha_2}
	   \centering
	   \includegraphics[width=\textwidth]{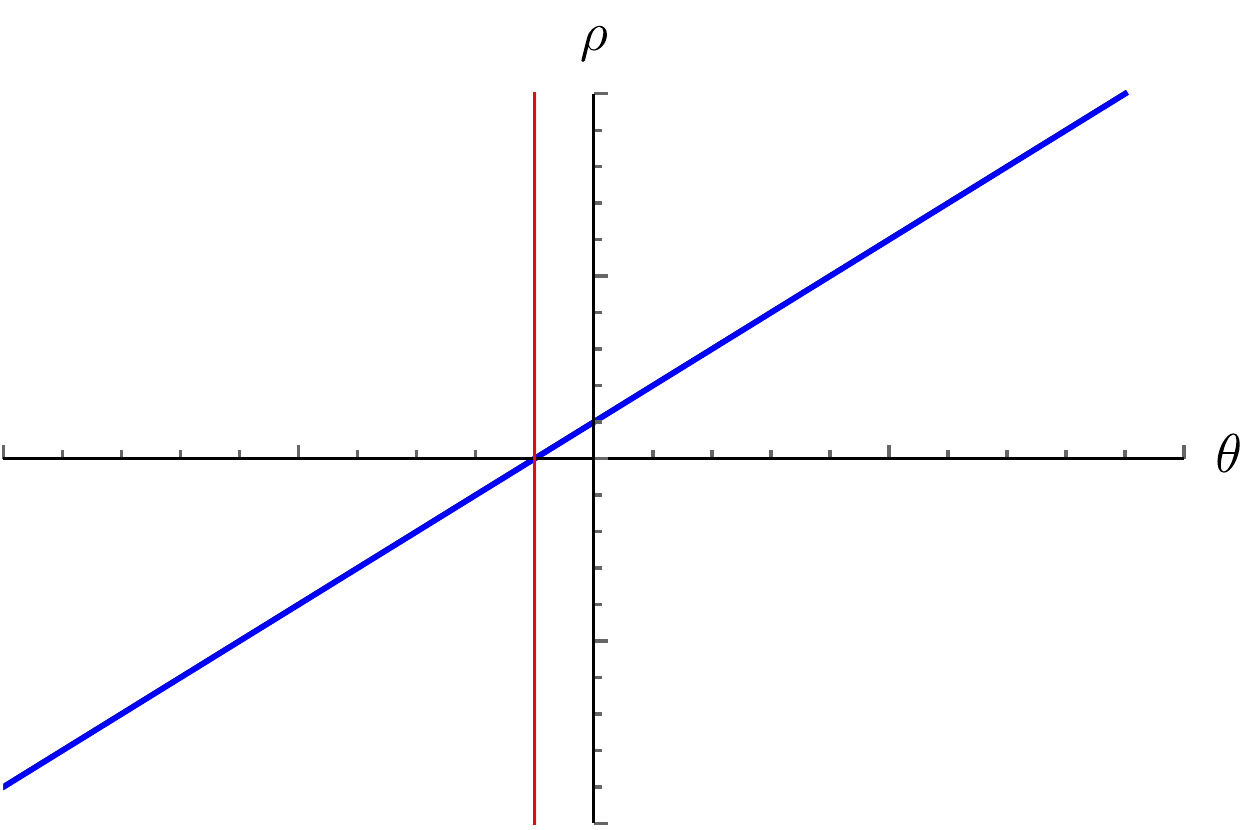}
	\end{minipage}}
\caption{Radius of curvature plots, without applying the bound of $\theta$, where the red line is the original bound. We prefer the case of Figure (b).}
\end{figure}

Since \eqref{eq:log-aesthetic_point_theta} includes $\rho$ of \eqref{eq:radius_curvature_theta}, let us define the extension using the radius of curvature that is $\left(\frac{\theta  \Lambda }{2}+1\right)^2$ in case of $\alpha=1.5$ (see \figref{fig:rho_alpha_1_5}), and $\theta  \Lambda +1$ in case of $\alpha=2$ (see \figref{fig:rho_alpha_2})

To resolve the ambiguity and define the log-aesthetic curve on the entire domain of $\theta$ similarly as in \figref{fig:log-aesthetic_alpha_2} when $\alpha > 1$, we mirror the radius of curvature across the line defined by the bound point (see red line \figref{fig:rho_alpha_2}) and also across the $\theta$-axis. Therefore, when $\alpha > 1$ the extended function of the radius of curvature is
\begin{equation} \label{eq:ext_radius_curvature_theta}
    \rho_{ext}^{\alpha>1}(\theta)=
\begin{cases}[c]
    \big((\alpha-1)\Lambda\theta+1\big)^{\frac{1}{\alpha-1}}
        & \theta > B_\theta\\
    -\big((1-\alpha ) \theta  \Lambda -1\big)^{\frac{1}{\alpha -1}}
        & \theta \leq B_\theta,
\end{cases}
\end{equation}
where $B_\theta=\frac{1}{\Lambda(1-\alpha)}$ is the earlier bound of $\theta$.
Based on this, a point of the extended log-aesthetic curve whose tangential angle is $\theta$ is defined on the complex plane as
\begin{equation} \label{eq:ext_log-aesthetic_theta}
    P_{ext}^{\alpha>1}(\theta) =
\begin{cases}[c] \nonumber
    \int_{0}^{\theta}  \big((\alpha-1) \Lambda \psi +1\big)^{\frac{1}{\alpha-1}} \me^{\iu\psi}\  d\psi
        & \theta>B_\theta \\
    \int_{0}^{\theta}  -\big((1-\alpha ) \Lambda \psi  -1\big)^{\frac{1}{\alpha -1}} \me^{\iu\psi}\  d\psi
        & \theta\leq B_\theta. 
\end{cases} 
\end{equation}
At $\theta=\frac{1}{\Lambda(1-\alpha)}$ the extended log-aesthetic curve includes a cusp.

In case of $\alpha=1$, the integration range of the algorithm is $[0,\theta_\Delta]$ and the log-aesthetic curve has no bound for $\theta$ nor $s$ hence there is no need to extend the log-aesthetic curve.

In case of $\alpha<1$, the integration range is also $[0,\theta_\Delta]$ and $\theta$ has an upper bound of $\frac{1}{1-\alpha}$. However, it does not yield problem until $\alpha \geq 0$ because at the point of $\theta=\frac{1}{1-\alpha}$ the arc length $s$ is infinite. The experimental result of \cite{yoshida2006interactive} also reports large drawable regions in this case. Therefore, there is no need to extend the log-aesthetic curve when $0 \leq \alpha < 1$.

\begin{figure}[!ht]
  \centering
  \subfloat[$\alpha=-0.5$]{
	\begin{minipage}[c][1\width]{0.45\textwidth} \label{fig:log-aesthetic_alpha_-0_5}
	   \centering
	   \includegraphics[width=\textwidth]{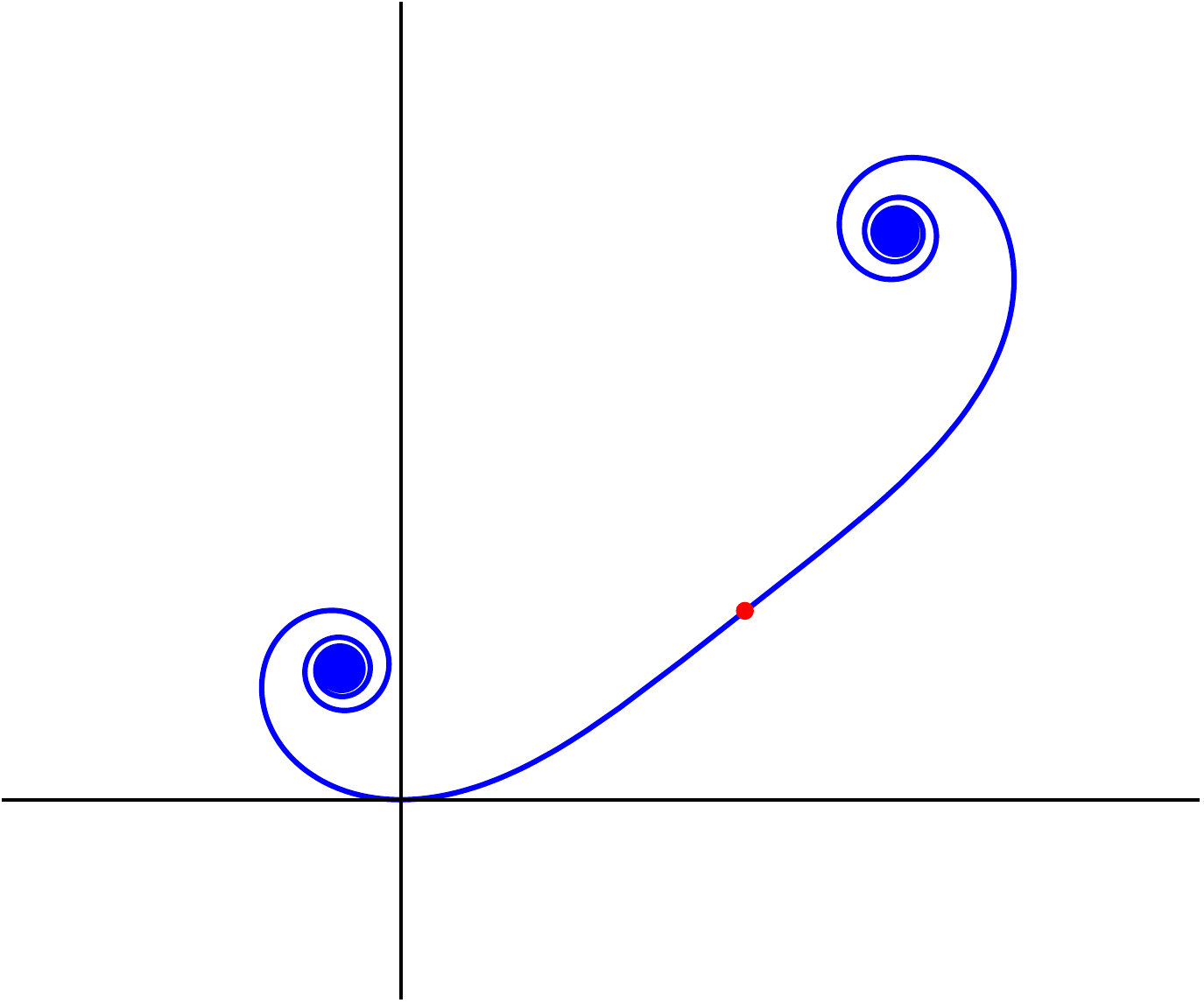}
	\end{minipage}}
 \hfill 	
  \subfloat[$\alpha=-1$]{
	\begin{minipage}[c][1\width]{0.45\textwidth} \label{fig:log-aesthetic_alpha_-1}
	   \centering
	   \includegraphics[width=\textwidth]{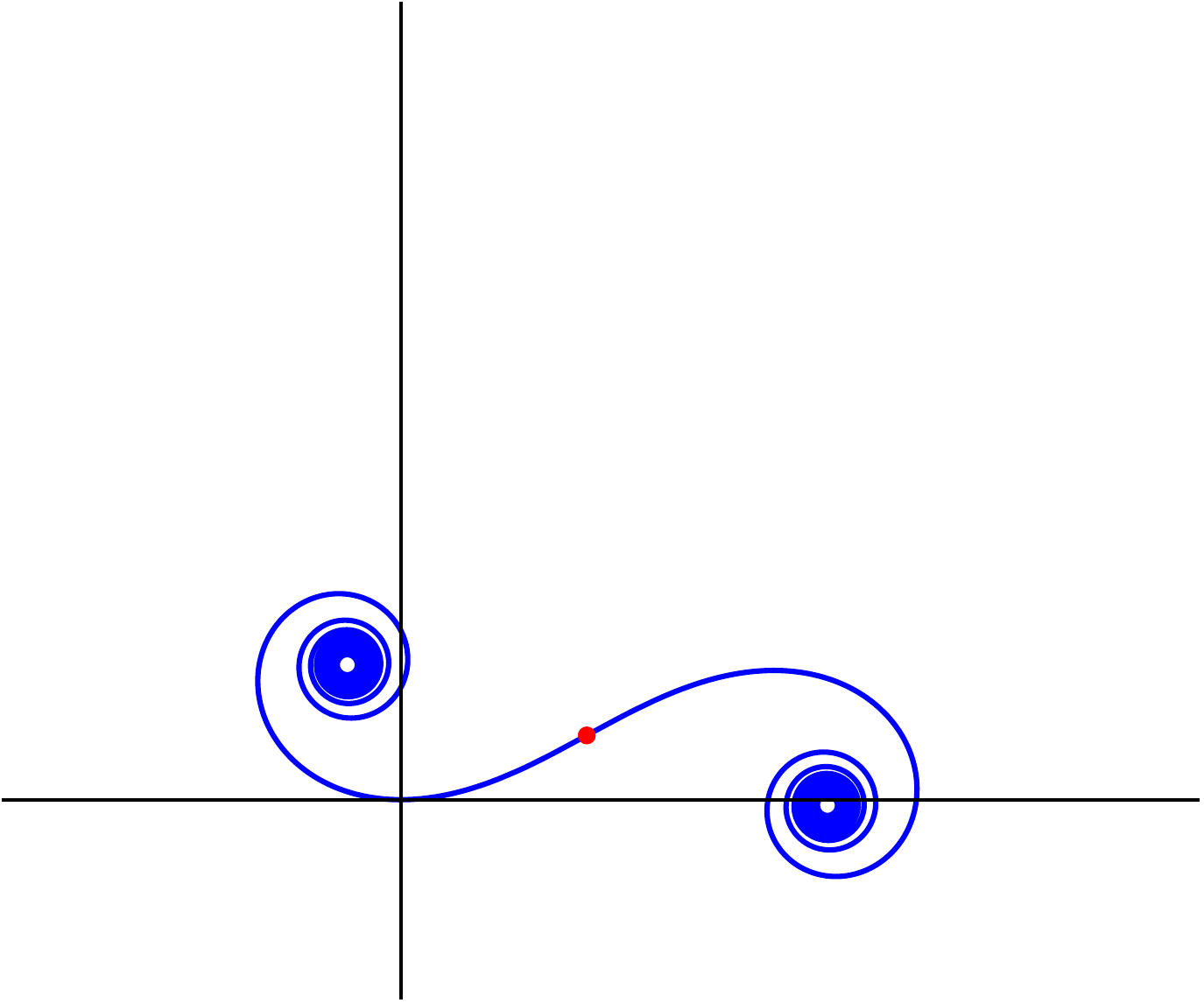}
	\end{minipage}}
\caption{Different log-aesthetic curves, without applying the bound of $s$.  We prefer the case of Figure (b).}
\end{figure}

\begin{figure}[!ht]
  \centering
  \subfloat[$\alpha=-0.5$]{
	\begin{minipage}[c][1\width]{0.45\textwidth}
	\label{fig:rho_alpha_-0_5}
	   \centering
	   \includegraphics[width=\textwidth]{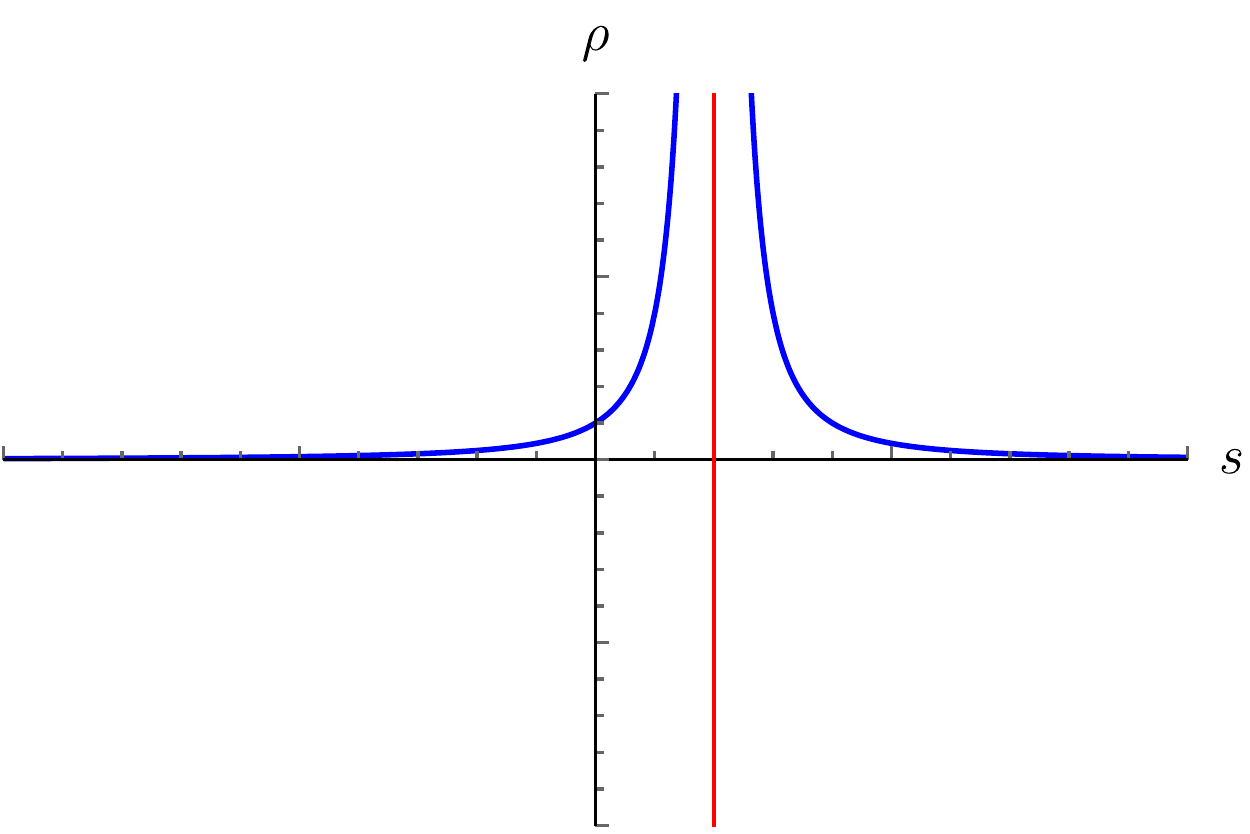}
	\end{minipage}}
 \hfill 	
  \subfloat[$\alpha=-1$]{
	\begin{minipage}[c][1\width]{0.45\textwidth}
	\label{fig:rho_alpha_-1}
	   \centering
	   \includegraphics[width=\textwidth]{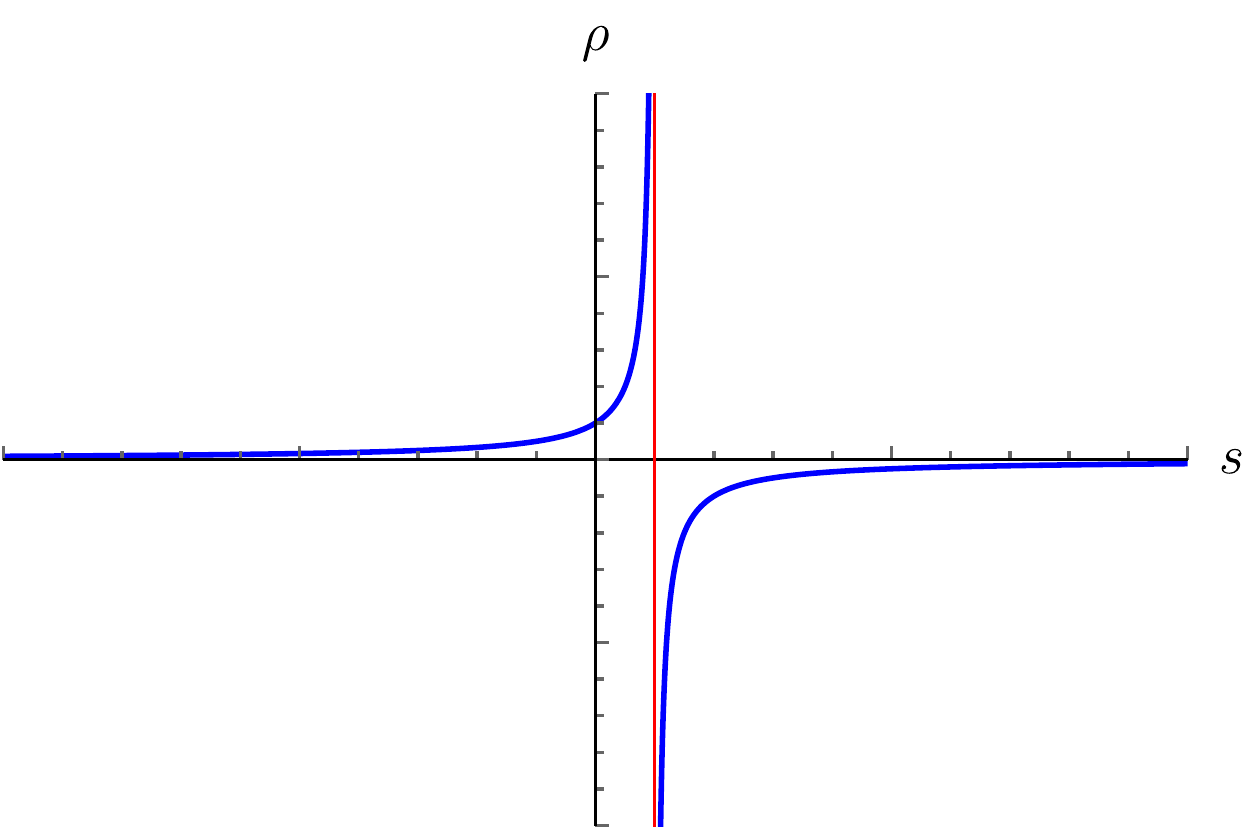}
	\end{minipage}}
\caption{Radius of curvature plots, without applying the bound of $s$, where the red line is the original bound. We prefer the case of Figure (b).}
\label{fig:rho_alpha_-0_5_-1}
\end{figure}

When $\alpha < 0$, the arc length $s$ also has an upper bound of $-\frac{1}{\Lambda \alpha}$. 
The integration range of the algorithm is $[0,\theta_\Delta]$ and the drawable region is limited by the point of bound. The extension of the log-aesthetic curve is necessary.
Regrading \cite{yoshida2006interactive}, the log-aesthetic curve has inflection point when $\alpha<0$ at the point of bound ($\theta=\frac{1}{\Lambda(1-\alpha)}$). Since $\rho$ is infinite (see \figref{fig:rho_alpha_-0_5_-1}) and the arc length is finite ($s=-\frac{1}{\alpha \Lambda}$), we can only use \eqref{eq:log-aesthetic_point_arc} to draw this point, and beyond we also apply the bound of \tableref{table:arc_theta_bounds}. However, there are also some cases depending on the value of $\alpha$ when the curve can be calculated on the entire domain of $s$.
For example, in case of $\alpha=-0.5$ (see \figref{fig:log-aesthetic_alpha_-0_5}), the \eqref{eq:log-aesthetic_point_arc} is
\begin{equation} \label{eq:log-aesthetic_arc_alpha_-0_5}
    C(s)= 
    \int_{0}^{s}\exp(\iu \frac{1}{12} \Lambda  u^2 (\Lambda  u-6)+u)\  d u,
\end{equation}
or in case of $\alpha=-1$ (when the log-aesthetic curve is a clothoid curve \cite{yoshida2006interactive}, see \figref{fig:log-aesthetic_alpha_-1}) it is
\begin{equation} \label{eq:log-aesthetic_arc_alpha_-1}
    C(s)= 
    \int_{0}^{s}\exp(\iu u-\frac{\Lambda  u^2}{2})\  d u.
\end{equation}

The different shapes depending on the value of $\alpha$ can be seen in the radius of curvature plots as well. When $\alpha=-0.5$, it is $\frac{4}{(\Lambda  s-2)^2}$ (see \figref{fig:rho_alpha_-0_5}) and it is $\frac{1}{1-\Lambda  s}$ in case of $\alpha=-1$ (see \figref{fig:rho_alpha_-1}).

We intend to extend the curve to increase the drawable region of the modeling algorithm, therefore, the curve is required to include an inflection point (as in case of $\alpha=-1$).
The reflection is defined based on the radius of curvature by mirroring its plot doubly across the line defined by the bound (see red line \figref{fig:rho_alpha_-1}) and across the $s$-axis.
Therefore, when $\alpha < 0$ the extended function of the radius of curvature is
\begin{equation} \label{eq:ext_radius_curvature_s}
    \rho_{ext}^{\alpha<0}(s)=
\begin{cases}[c]
    (\alpha  \Lambda  s+1)^{\frac{1}{\alpha }}
        & s < B_s\\
    -(-\alpha  \Lambda s -1)^{\frac{1}{\alpha }}
        & s \geq B_s,
\end{cases}
\end{equation}
where $B_s=-\frac{1}{\alpha  \Lambda }$ is the earlier bound of $s$.
A point of the extended log-aesthetic curve whose arc length is $s$ defined as:
\begin{equation} \label{eq:ext_log-aesthetic_arc}
    C_{ext}^{\alpha<0}(s)= \\
\begin{cases}[c]
    \int_{0}^{s}\exp(\iu \frac{(\alpha  \Lambda  u+1)^{\frac{\alpha -1}{\alpha }}-1}{(\alpha -1) \Lambda })\  d u
        & s \leq B_s \\
    \int_{0}^{s}\exp(\iu \frac{(\alpha  \Lambda  (-u)-1)^{1-\frac{1}{\alpha }}-1}{(\alpha -1) \Lambda })\  d u
        & s > B_s,
\end{cases}
\end{equation}
where the extended log-aesthetic curve includes an inflection point.

The extensions of the log-aesthetic curves in the above cases describe two different mirrorings (one by tangential angle and the other by arc length), that defines two different equations, which need to be used: \eqref{eq:ext_log-aesthetic_theta} when $\alpha > 1$ and \eqref{eq:ext_log-aesthetic_arc} when $\alpha < 0$.
Otherwise, we need to apply the bound of $\theta$ and $s$ of \tableref{table:arc_theta_bounds}. Namely, there is an upper bound of $\theta$ and a lower bound of $s$. Therefore, we use \eqref{eq:log-aesthetic_point_theta} (the original formula by tangential angle \cite{yoshida2006interactive}) in case of $s \leq 0$ (and $\theta \leq 0$ as well) and we use \eqref{eq:log-aesthetic_point_arc} (original equation by arc length \cite{yoshida2006interactive}) in case of $s > 0$ (and $\theta > 0$), when $0 \leq \alpha \leq 1$.

However, the equations by arc length require to determine $s$ from the given $\theta_\Delta$. To do it when $0 \leq \alpha \leq 1$, the following equation can be used (from \eqref{eq:theta_arc}):
\begin{equation} \label{eq:arc_theta}
    \theta(s)^{-1} = S(\theta) = 
\begin{cases}[c]
    -\frac{\log (1-\theta  \Lambda )}{\Lambda } 
        & \text{if } \alpha = 0\\
    \frac{e^{\theta  \Lambda }-1}{\Lambda }
        & \text{if } \alpha = 1\\
    \frac{((\alpha -1) \theta  \Lambda +1)^{\frac{\alpha }{\alpha -1}}-1}{\alpha  \Lambda }  
         & \text{if } 0 < \alpha < 1. \\
\end{cases}
\end{equation}

In case of $\alpha < 0$, the derivation of the formula is not straightforward.
Although, in \eqref{eq:arc_theta} $\theta$ has upper bound (except when $\alpha=1$) the arc length $s$ increases up to infinity ($S(\frac{1}{\Lambda(1-\alpha)})=\infty$).
The arc length is extended when $\alpha<0$ (\eqref{eq:ext_log-aesthetic_arc}), so it has no upper bound as in the original case, however, the tangential angle is increasing before the inflection point and decreasing after (see \figref{fig:log-aesthetic_alpha_-1}). That means, $\theta$ still has an upper bound. Moreover, the $S(\theta)$ is not a one to one correspondence in case of $\alpha<0$. Therefore it can be defined either as
\begin{equation} \label{eq:arc_theta_within}
    S(\theta)_{within}^{\alpha<0} = 
        \frac{((\alpha -1) \theta  \Lambda +1)^{\frac{\alpha }{\alpha -1}}-1}{\alpha  \Lambda }
\end{equation}
or
\begin{equation} \label{eq:arc_theta_beyond}
    S(\theta)_{beyond}^{\alpha<0} = 
        -\frac{((\alpha -1) \theta  \Lambda +1)^{\frac{\alpha }{\alpha -1}}+1}{\alpha  \Lambda },
\end{equation}
depending on whether the point of the curve is within the inflection point or beyond it. Since the log-aesthetic curves are mirrored parts, the following is satisfied:
\begin{equation} \label{eq:arc_theta_equality}
    S(\theta)_{within}^{\alpha<0} = S(2 B_\theta - \theta)_{beyond}^{\alpha<0},
\end{equation}
where $B_\theta=\frac{1}{(1-\alpha)\Lambda}$ (the upper bound of $\theta$) is the point of inflection.

\subsection{Finding the Curve Segment Using a New \texorpdfstring{\boldmath$\Lambda$}{Lambda} Bisection Method}
\label{sec:lambda_bisection}

In the new approach, a similar $\Lambda$ bisection method is used as in \cite{yoshida2006interactive} to find the extended log-aesthetic curve segment that fits the triangle $ABC$. However, in the new algorithm, the desired curve segment is specified by two points ($A$ and $C$) and two vectors $\vec{v_A}$ and $\vec{v_C}$, where $\vec{v_A}$ is the tangent vector of the curve at $A$, and $\vec{v_C}$ defines only the direction of the tangent line at $C$.
The point $B$ is the intersection point of the tangent lines. The difference of the tangential angle between the first and last endpoint ($\theta_\Delta$) is obtained by calculating the angle $\beta$ between $\vec{v_A}$ and $\vec{v_C}$, thus $\theta_\Delta = \pi - \beta$. 

The bisection method (as well as in \cite{yoshida2006interactive}) repeatedly bisects the interval defined for $\Lambda$ and selects the sub-interval in which the absolute difference between the compared angles of the two triangles becomes smaller.
In the algorithm of \cite{yoshida2006interactive}, in case of $\alpha \leq 1$, the angle $\theta_A$ is compared with $\theta_{A'}$, and the angle $\theta_C$ is compared with $\theta_{C'}$ when $\alpha > 1$.
Besides, the coordinates of point $A$ and $C$ are swapped if $\left| A B \right| \leq \left| B C \right|$ does not hold.
In the new algorithm, a flag is marked to indicate the endpoint swapping and the bisection algorithm decides which angle to use ($\theta_A$ or $\theta_C$)  to compare the triangles.

Since the log-aesthetic curve is not extended when $0 \leq \alpha \leq 1$, the original method is used without change.
Otherwise, when $\alpha < 0$ or $\alpha > 1$ the bisection method requires several modifications.

In case of $\alpha<0$, the tangential angle increases from the reference point until the inflection point ($\theta=\frac{1}{\Lambda(1-\alpha)}$) and decreases beyond it. This means that the interval for $\Lambda$ is still between $0$ and $1/(\theta_\Delta(1-\alpha)$, but we need to test whether point $C'$ is within the bound or beyond it, since the selection of the sub-intervals of the $\Lambda$ bisection needs to be changed conversely. At this point, the algorithm decides which equation to use to determine arc length $s$ from $\theta_\Delta$, \eqref{eq:arc_theta_within} or \eqref{eq:arc_theta_beyond}.

On the other hand, when $\alpha > 1$, the tangential angle decreases from the reference point ($\theta=0$) to the singular point ($\theta=\frac{1}{\Lambda(1-\alpha)}$) and still reducing beyond it. Therefore, the extension of the log-aesthetic curve means also the extension of the $\Lambda$ range during the bisection to go beyond the cusp with $A'$. Thus, the new interval is $0 < \Lambda < 2/(\theta_\Delta(1-\alpha))$ in this case.

\begin{figure}[ht]
  \subfloat[$\Lambda = 0.588$]{
	\begin{minipage}[t]{0.49\textwidth}
	   \centering
	   \includegraphics[width=0.9\textwidth]{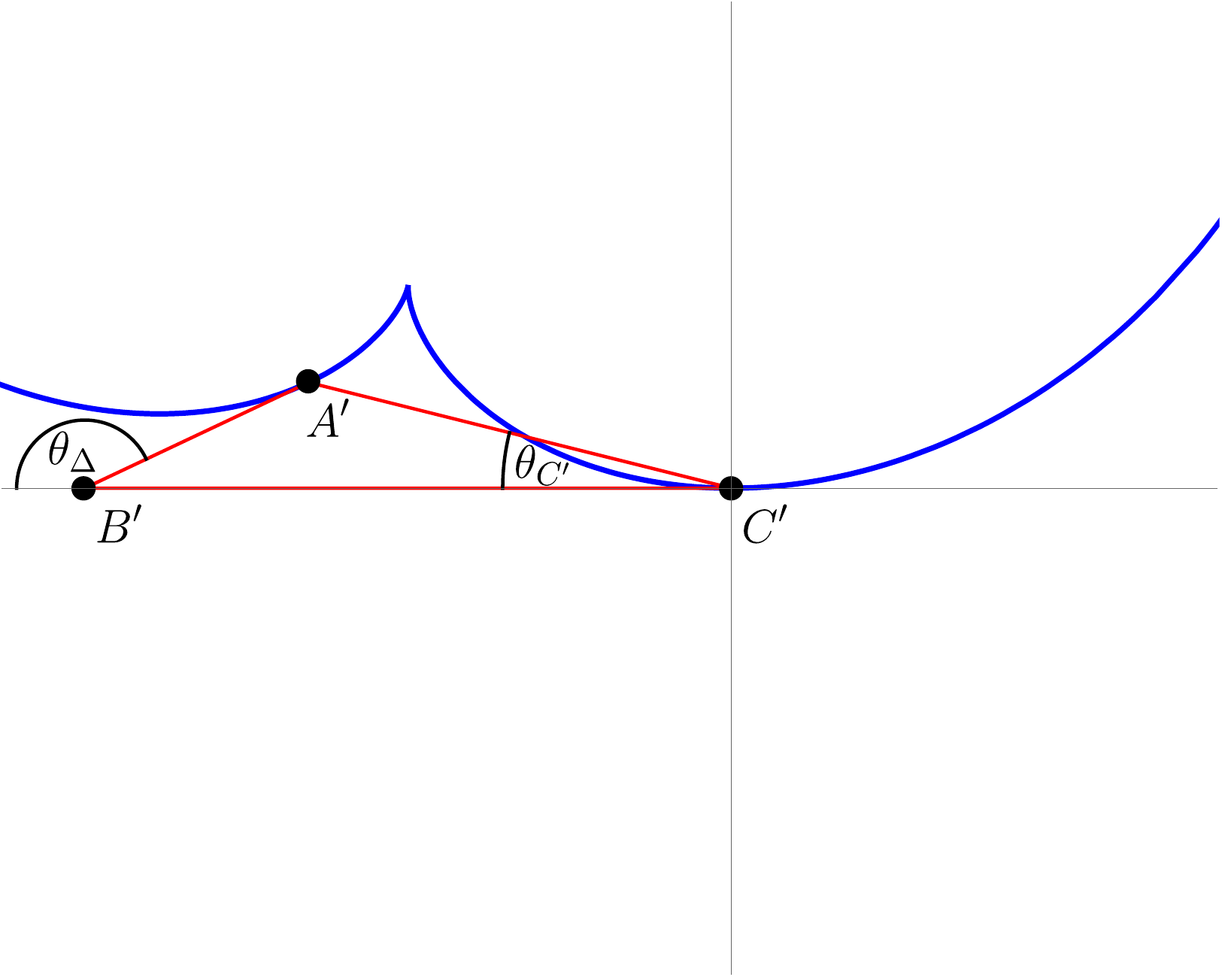}
	\end{minipage}}
 \hfill 	
  \subfloat[$\Lambda = 0.87$]{
	\begin{minipage}[t]{0.49\textwidth}
	   \centering
	   \includegraphics[width=0.9\textwidth]{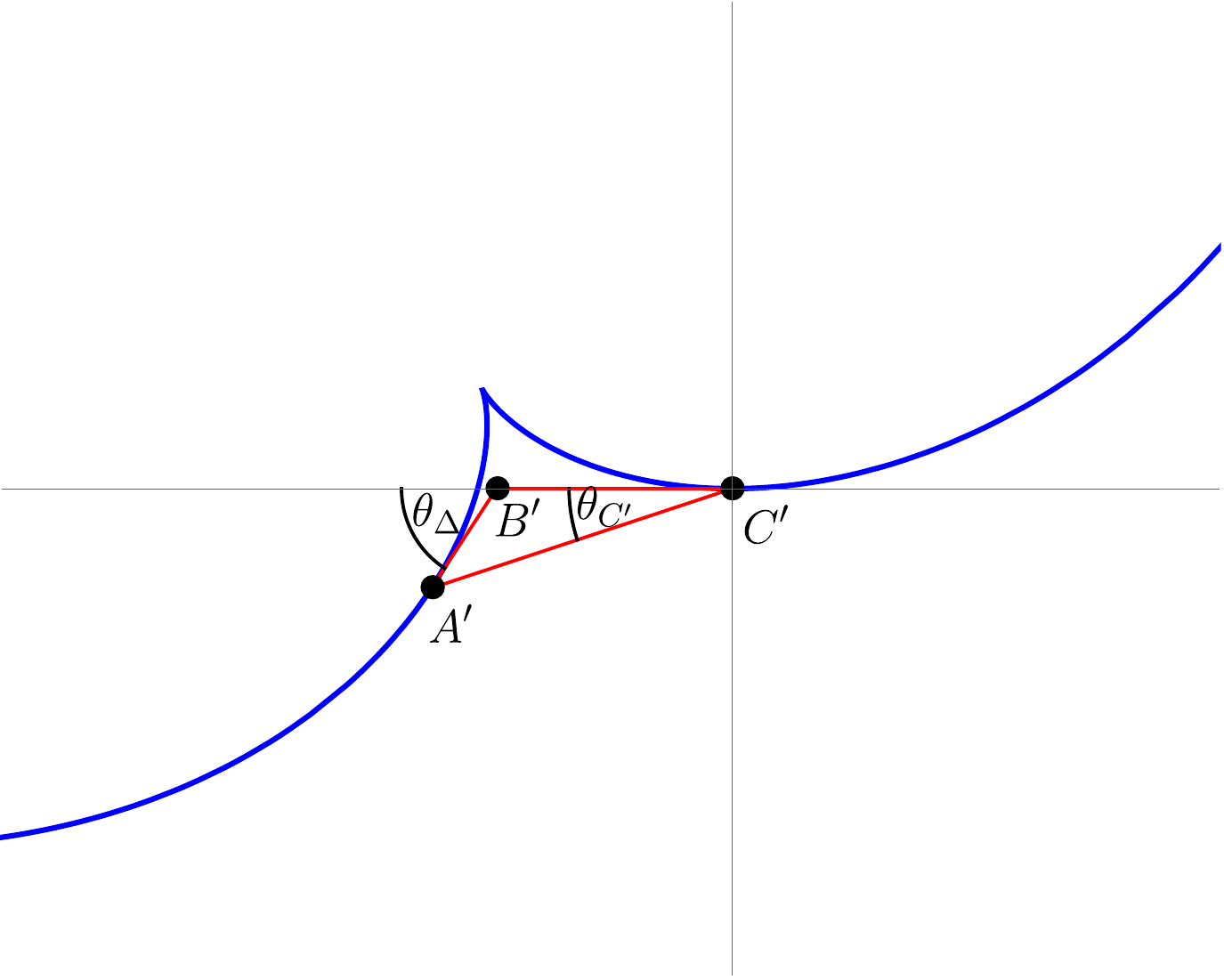}
	\end{minipage}}
\caption{Different orientations of the triangle $A'B'C'$ during the $\Lambda$ bisection method depending on the position of $A'$ ($\alpha=2$).}
\label{fig:different_orientation_plus}
\end{figure}

In the original algorithm of \cite{yoshida2006interactive}, in case of $\alpha > 1$, when $\Lambda$ is increasing the $\theta_{C'}$ is also increasing. However, if point $A'$ is below the \mbox{$x$-axis} (that can happen with the extended log-aesthetic curves) the $\theta_{C'}$ is decreasing because the orientation of the triangle $A'B'C'$ is changed. Therefore, the appropriate sub-interval need to be selected depending on the position of point $A'$ (see \figref{fig:different_orientation_plus}).

\begin{figure}[!ht]
  \subfloat[$\Lambda = 0.145$]{
	\begin{minipage}[t]{0.49\textwidth}
	   \centering
	   \includegraphics[width=0.9\textwidth]{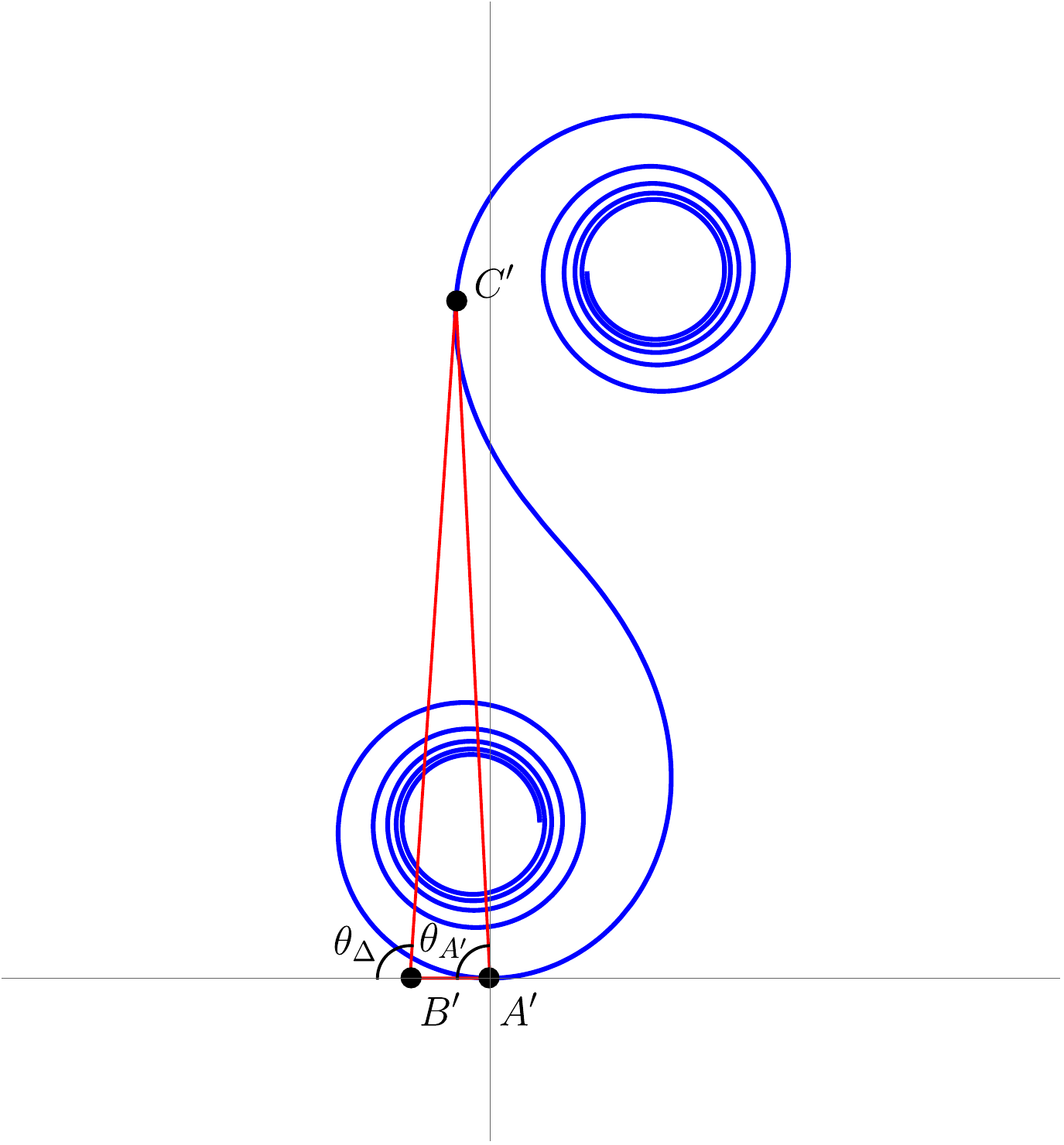}
	\end{minipage}}
 \hfill 	
  \subfloat[$\Lambda = 0.17$]{
	\begin{minipage}[t]{0.49\textwidth}
	   \centering
	   \includegraphics[width=0.9\textwidth]{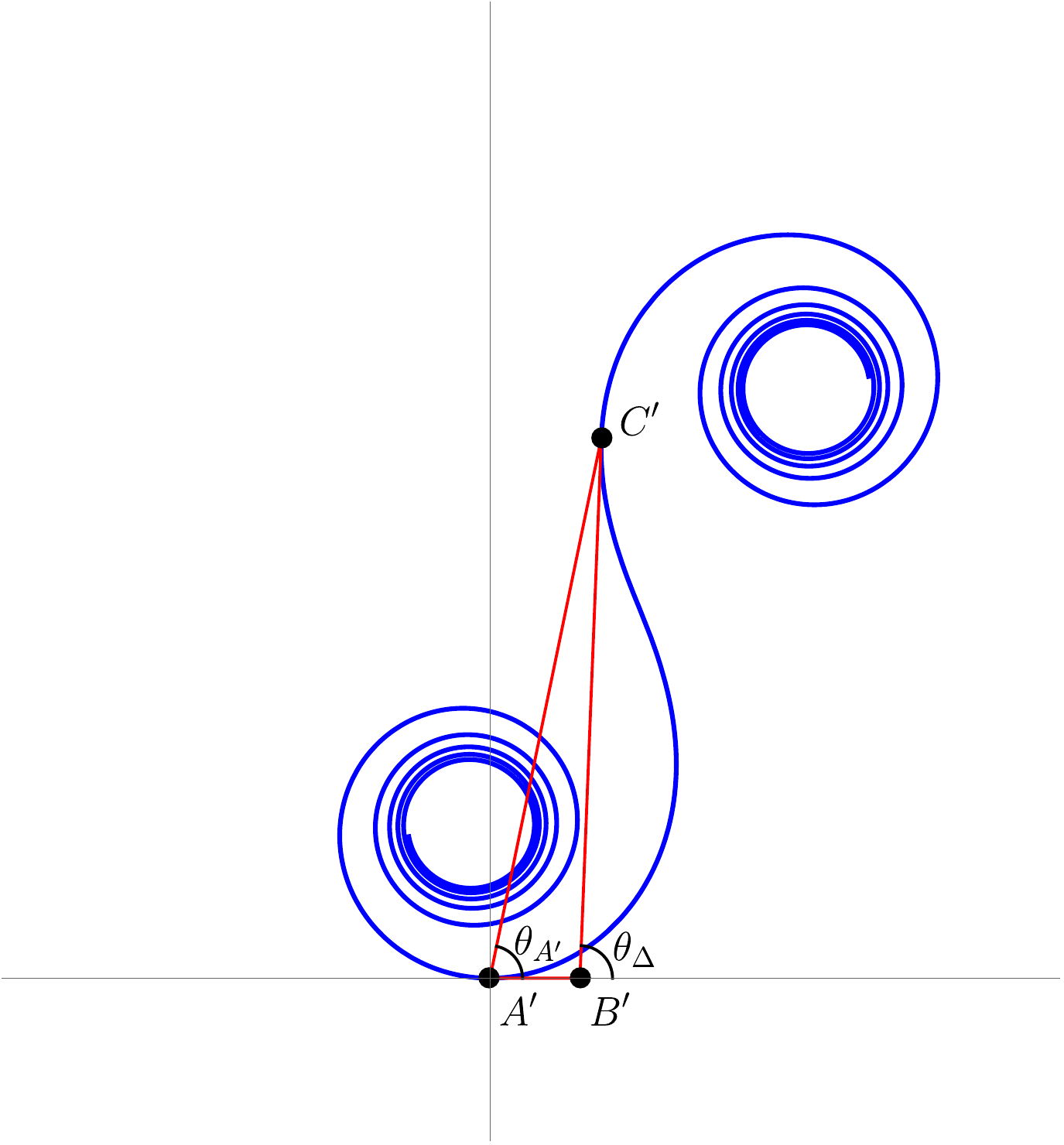}
	\end{minipage}}
\caption{Different orientations of the triangle $A'B'C'$ during the $\Lambda$ bisection method depending on the position of $B'$ ($\alpha=-2$).}
\label{fig:different_orientation_minus}
\end{figure}

In case of $\alpha < 0$, besides the above case, another event may happen. In the original algorithm of \cite{yoshida2006interactive}, the tangential lines of $A'$ and $C'$ always cross each other on the right side of the \mbox{$y$-axis}. In the new $\Lambda$ bisection method, the cross point $B'$ can also appear on the left side of the \mbox{$y$-axis} because of the extension, and it also changes the orientation of the triangle $A'B'C'$ (see \figref{fig:different_orientation_minus}).

For the complete pseudo-code of the new $\Lambda$ bisection algorithm that implements the extended range and solves the orientation problems see \listref{lst:lambda_bisection}.
The modifications regarding the previous algorithm of \cite{yoshida2006interactive} are highlighted in red.
The presented algorithm is able to draw the log-aesthetic curve with only a minor condition, considering the case when $0 \leq \alpha < 1$. However, regarding both \cite{yoshida2006interactive} and \cite{yoshida2012drawable}, the small restriction for the placement is possibly caused by the computation error of the large integration range.

The authors of \cite{yoshida2012drawable} are also presented an algorithm to draw the theoretical drawable region of log-aesthetic curve segments.
However, these regions are expanded with the new $\Lambda$ bisection method and using the extended log-aesthetic curves, their algorithm still can be used during the modelling process to indicate the possible location of the inflection point (when $\alpha < 0$) or the curvature-extremal point (when $\alpha > 1$) depending on the position of the given geometric data.

\subsection{Controlling the Tangent Length and Curvature of the First Point Using an \texorpdfstring{\boldmath$\alpha$}{Alpha} Bisection Method}
\label{sec:alpha_bisection}
    
The presented $\Lambda$ bisection with the given geometric data determines the extended log-aesthetic curve segment with an arbitrary value of $\alpha \in \R$. Regarding Harada et al. \cite{harada1999aesthetic}, $\alpha$ is related to the impression of the curve. However, it is difficult to choose a suitable value for it and it is a common practice to fix the parameter to design with log-aesthetic curves.
In this subsection, we present an algorithm to determine $\alpha$ to match also the length of $\vec{v_A}$ and control the radius of curvature at the first point with the length of its tangent.
Since the same geometric data with different $\alpha$ parameters require different $\Lambda$ values, an exact calculation is not possible. Therefore, we use another bisection method to determine the appropriate value of $\alpha$.

\begin{figure}[ht]
  \subfloat[$\alpha \approx -3.4, \Lambda \approx 0.2$, \lstinline{swap_flag} = false]{
	\begin{minipage}[t]{0.49\textwidth}
	\label{fig:final_alpha_minus_not_swapped}
	   \centering
	   \includegraphics[width=0.95\textwidth]{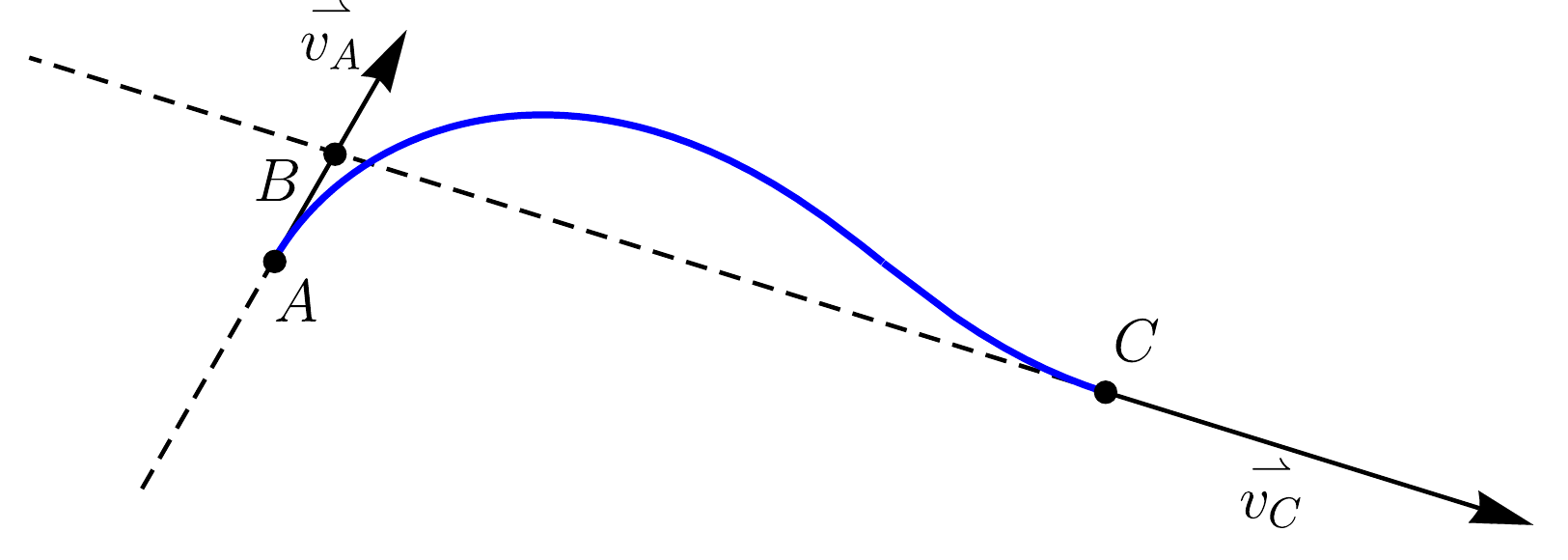}
	\end{minipage}}
 \hfill 	
  \subfloat[$\alpha \approx 1.66, \Lambda \approx 1.12$, \lstinline{swap_flag} = false]{
	\begin{minipage}[t]{0.49\textwidth}
	\label{fig:final_alpha_plus_not_swapped}
	   \centering
	   \includegraphics[width=0.95\textwidth]{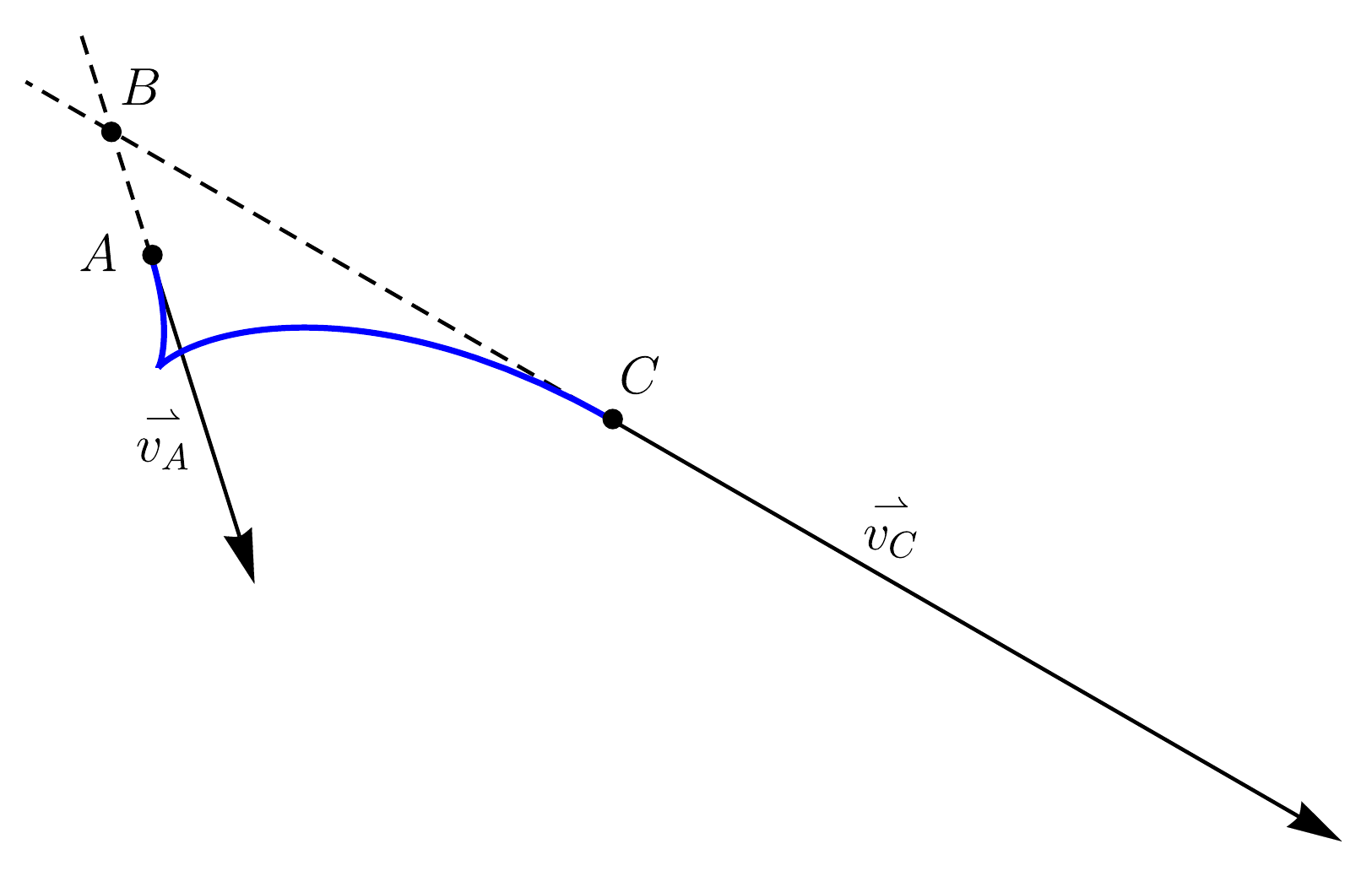}
	\end{minipage}}\\
  \subfloat[$\alpha \approx -0.8, \Lambda \approx 0.3$, \lstinline{swap_flag} = true]{
	\begin{minipage}[t]{0.49\textwidth}
	\label{fig:final_alpha_minus_swapped}
	   \centering
	   \includegraphics[width=0.95\textwidth]{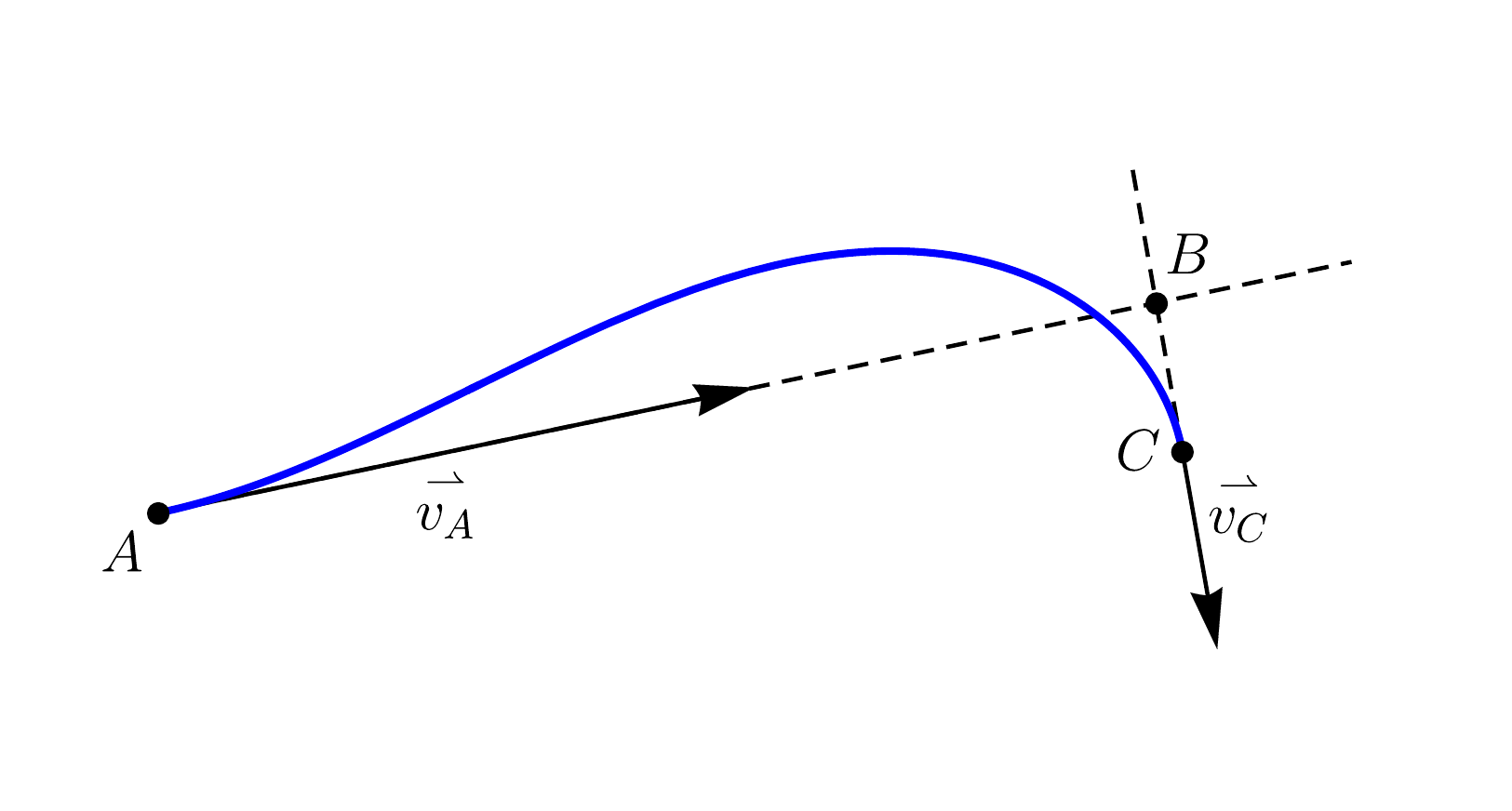}
	\end{minipage}}
 \hfill 	
  \subfloat[$\alpha \approx 2.5, \Lambda \approx 0.6$,
  \lstinline{swap_flag} = true]{
	\begin{minipage}[t]{0.49\textwidth}
	\label{fig:final_alpha_plus_swapped}
	   \centering
	   \includegraphics[width=0.95\textwidth]{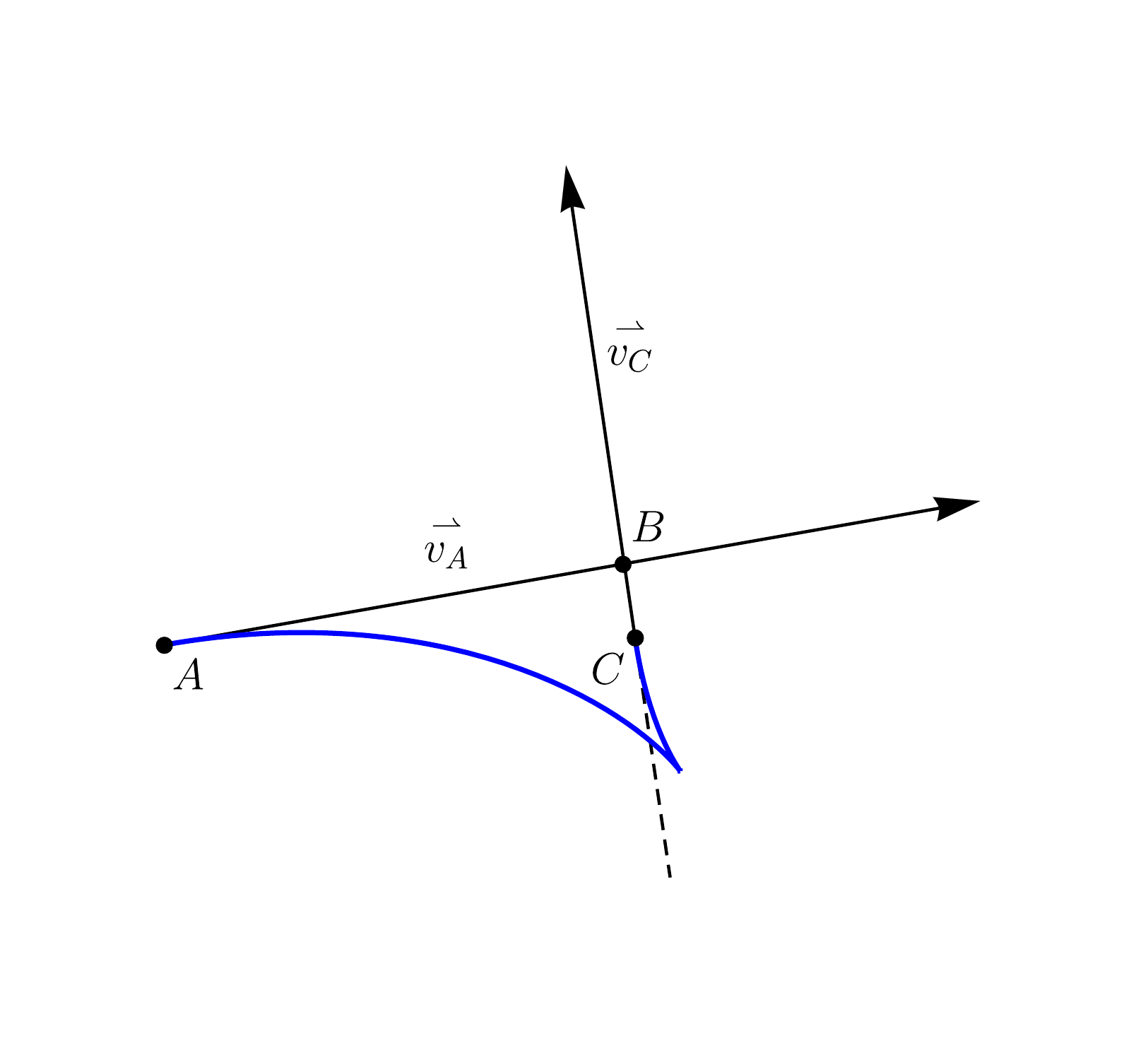}
	\end{minipage}}
\caption{Example results of the $\alpha$ bisection. The drawn curve segments are defined by the given points $A$ and $C$, the vector $\vec{v_A}$, and the direction of $\vec{v_C}$.}
\label{fig:final_example}
\end{figure}

There are two different instances whether point $B$ is closer to $A$ or $C$, the \lstinline{swap_flag} is false or true.
In the first case, when the \lstinline{swap_flag} is false the extended log-aesthetic curve segment may include inflection point if $\vec{v_A}$ points to $B$ from $A$, and it may have cusp otherwise (see \figref{fig:final_alpha_minus_not_swapped} and \figref{fig:final_alpha_plus_not_swapped}).
On the other hand, when the \lstinline{swap_flag} is true, the determined log-aesthetic curve segment is generated from $C$ to $A$. However, we desire point $A$ to be the first, therefore, to define the appropriate tangent directions it is preferred to apply a reversed parameter transformation on the log-aesthetic curve segment to obtain oppositely directed tangent vectors, as it is in \figref{fig:final_alpha_minus_swapped} and \figref{fig:final_alpha_plus_swapped}.In this case, if $\vec{v_C}$ points to $B$ from $C$, the curve segment may have cusp and it can include inflection point otherwise.

Let us see the cases when the coordinate \lstinline{swap_flag} is false, the cross point $B$ is closer to $A$ ($\left| A B \right| \leq \left|C B\right|$).
When $\alpha \leq 1$, $A'$ is the reference point, where the log-aesthetic curve has unit tangent vector (see \figref{fig:previous_algorithm_conf_1}). Therefore, the tangent length of $A'$ is always unit after the $\Lambda$ bisection. It corresponds to the given $A$, where the tangent vector is $\vec{v_A}$.
Thus, its length equals the scale factor of the geometric transformation (that transforms the log-aesthetic points from $A'B'C'$ to $ABC$).
The scale value can be computed as the ratio of the endpoint distances: $\left|A C\right|/\left|A'C'\right|$. When $\alpha$ decreases, the points $A'$ and $C'$, under the same $\theta_\Delta$ become closer on the plane. It is because the log-aesthetic curve in case of $\alpha < 1$ spirally converges to the point where $\rho=0$ as $\theta$ approaching $-\infty$. The convergence is faster on lower values of $\alpha$.
As a result, the minimum of $\left|A'C'\right|$ is when $\alpha=-\infty$.
Therefore, the scale factor of the transformation and consequently the length of vector $\vec{v_A}$ is the highest in this case.
 
When $\alpha > 1$ (and $\left| A B \right| \leq \left|C B\right|$ still holds), the given $A$ does not correspond to the reference point, $A'$ is defined as the log-aesthetic point where the tangential angle is $-\theta_\Delta$ (see \figref{fig:previous_algorithm_conf_2}). Therefore, the tangent vector of this point has to be computed as the derivation of \eqref{eq:ext_log-aesthetic_theta}. It means that the length of the vector $\vec{v_A}$ does not depend on the scale factor only, but on the radius of curvature as well.
When $\alpha > 1$, the integration range is $[-\theta_\Delta,0]$, where $\rho$ (and the tangent length as well) decreases from the reference point. At the original bound of $\theta=\frac{1}{\Lambda(1-\alpha)}$ the $\rho=0$ and the tangent becomes zero. Since the bisection range of the $\Lambda$ is extended, the point can go beyond the cusp as well, where the sign of $\rho$ (and the direction of the tangent vector) is changed and the tangent length increases until it becomes unit again (at $\theta= \frac{2}{\Lambda(1-\alpha)}$). It means that when $\alpha > 1$ after the $\Lambda$ bisection in the range of $0 < \Lambda < 2/(\theta_\Delta(1-\alpha))$ the tangent length of $A'$ is in the interval of [-1,1]. Besides, this tangent length is multiplied by the scale factor of the geometric transformation, similarly as in the previous case. However, in case of $\alpha > 1$, the log-aesthetic curve spirally diverges towards the point where $\rho=\infty$. Thus, by enlarging the value of $\alpha$, the velocity of the divergence increases until $\alpha=\infty$ when the tangent length became unit again.

\begin{figure}[ht]
  \subfloat[$\alpha=-\infty$]{
	\begin{minipage}[t]{0.49\textwidth}
	\label{fig:alpha_bisect_not_swapped_max_plus}
	   \centering
	   \includegraphics[width=0.95\textwidth]{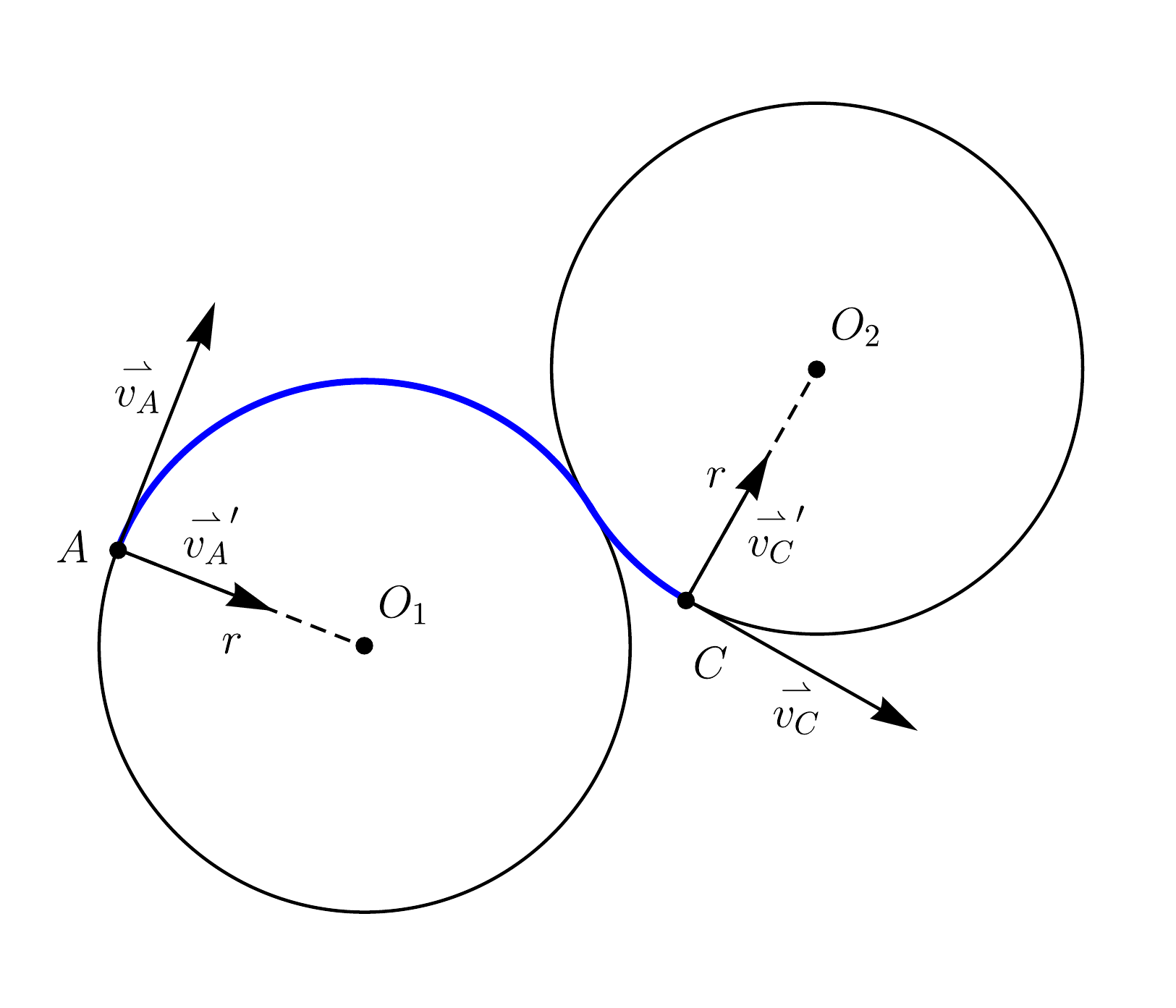}
	\end{minipage}}
 \hfill 	
  \subfloat[$\alpha=+\infty$]{
	\begin{minipage}[t]{0.49\textwidth}
	\label{fig:alpha_bisect_not_swapped_max_minus}
	   \centering
	   \includegraphics[width=0.95\textwidth]{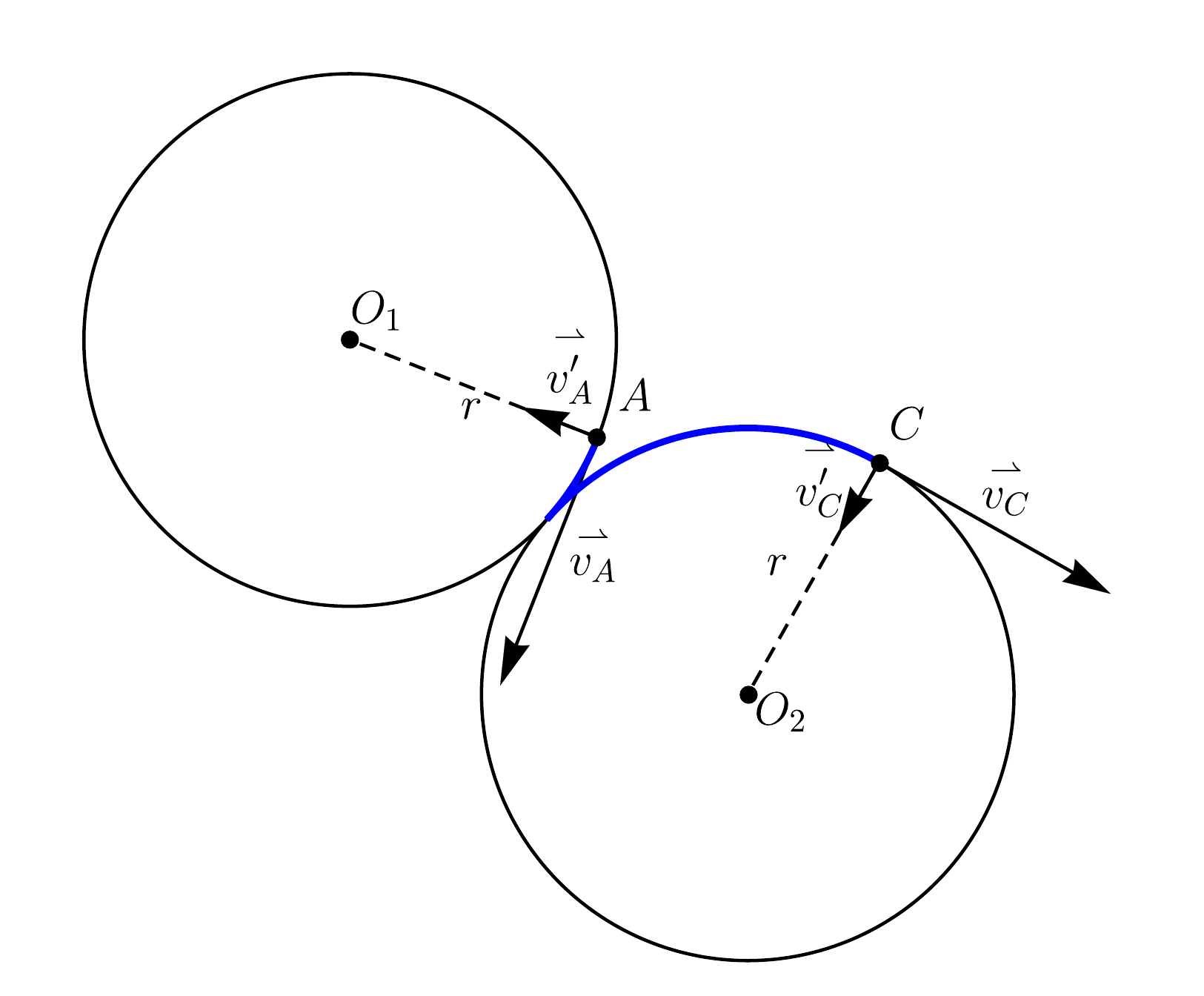}
	\end{minipage}}
\caption{Calculation of the minimum and maximum lengths of $\vec{v_A}$ and $\vec{v_C}$.}
\label{fig:alpha_bisect_not_swapped_max}
\end{figure}

It can be concluded, when the endpoints are not swapped there is a well-defined range for the length of the tangent vector at $A$ and the monotonic change of $\alpha$ varies the tangent simultaneously. At a certain value of $\alpha$, point $A$ becomes the cusp and the tangent vector $\vec{v_A}$ changes its direction. Below this $\alpha$ value the extended log-aesthetic curve segment may include inflection point and above may include cusp.
By calculating the superior and inferior length values, a new bisection method can be used on $\alpha$ to find its appropriate value to match the given length of $\vec{v_A}$, if it is between the range.
The maximum length is computed when $\alpha=-\infty$, when the log-aesthetic curve is two touching unit circles centered at $[0,1]$ and $[0,-1]$. In this case, the given points $A$ and $C$ are placed on two touching circles with unknown but equal radius $r$ that is the scale factor of the geometric transformation (See \figref{fig:alpha_bisect_not_swapped_max_plus}). The length of $\vec{v_A}$ and $\vec{v_C}$ is also $r$.
Let $O_1$ and $O_2$ be the center of these circles. Using the parametric equations of the line defined by the points $A$ and $O_1$ and $C$ and $O_2$, the following can be written:
\begin{align} \label{eq:computing_circles_radius}
    O_1 = A + r \cdot \vec{v_A}' \nonumber \\
    O_2 = C + r \cdot \vec{v_C}' \\
    \left| O_2 - O_1 \right| = 2 r, \nonumber
\end{align}
where $\vec{v_A}'$ and $\vec{v_C}'$ are unit vectors, directed along $\overline{AO_1}$ and $\overline{CO_2}$ (See \figref{fig:alpha_bisect_not_swapped_max_plus}).

Since $A$ and $C$ is given and $\vec{v_A}'$ and $\vec{v_C}'$ can be calculated, the solution of the equation system above gives the scale factor in case of $\alpha=-\infty$, the maximum length of the tangent vector at $A$.
When $\alpha=-\infty$ the extended log-aesthetic curve is two, fix centered circle, independently of the value of $\Lambda$. This means that the $\Lambda$ bisection has solution only when $\alpha > -\infty$. Therefore, the above value $r$ is the theoretical maximum length of the vector $\vec{v_A}$ that can never be reached.

On the other hand, the maximum length of $\vec{v_A}$ in the opposite direction is calculated when $\alpha=\infty$. In this case, the extended log-aesthetic curve become two touching circles once again. Therefore, a similar equation is used as \eqref{eq:computing_circles_radius}, however, the rotation direction of vector $\vec{v_A}$ and $\vec{v_C}$ need to be changed (see \figref{fig:alpha_bisect_not_swapped_max_minus}). The solution of this equation is also a theoretical maximum length of $\vec{v_A}$. 

In the previous cases, when the coordinate \lstinline{swap_flag} is false, there are a certain $\alpha$ ($>1$), when the first point is the cusp. Above and below this value, the tangent at $A'$ (and $A$ as well) has a different direction.
Therefore, the vector $\vec{v_A}$ is unique on different $\alpha$ values.
However, when $\left| A B \right| > \left| C B\right|$ and the \lstinline{swap_flag} is true, there is a certain $\alpha$ ($<1$) value, when the first point $A$ become the point of the inflection and the length of $\vec{v_A}$ is infinite. Above and below this value, the length is only decreasing but the vector has the same direction.
Therefore, there are two similar instances, in both the maximum length of the tangent vector $\vec{v_A}$ is infinite, when $A$ is the inflection point.
The minimum length is different depending on whether the $\alpha$ approaches $-\infty$ or $+\infty$. The vector $\vec{v_A}$ has the same direction in the two different cases and the extended log-aesthetic curve segment may include inflection or cusp otherwise. The minimum lengths of $\vec{v_A}$ are also calculated using \eqref{eq:computing_circles_radius} (See also \figref{fig:alpha_bisect_not_swapped_max_plus} and \figref{fig:alpha_bisect_not_swapped_max_minus}).

Since the bisection determines the value of $\alpha$, the algorithm decides which case to use depending on the given direction of $\vec{v_C}$. The point $C$ is on the opposite side of the cusp when $\alpha$ approaching $+\infty$ or $-\infty$ and its tangent vector $\vec{v_C}$ has opposite direction in the two different cases (see \figref{fig:final_alpha_minus_swapped} and \figref{fig:final_alpha_plus_swapped}). Therefore, if it points to $B$, the extended log-aesthetic segment may include cusp, and it may include inflection point otherwise.

During the $\alpha$ bisection, the unnecessary sub-intervals of $\alpha$ are excluded based on the position of inflection point, using \lstinline{beyond_inf_point} of \listref{lst:lambda_bisection}. Moreover, in the case when the inflection point is included, the sub-intervals needs to be selected contrary, since $\alpha$ needs to be decreased to approach the inflection point with $A$ and gain larger vector length at $A$.

For the pseudo-code of the $\alpha$ bisection algorithm see \listref{lst:alpha_bisection}. The function takes the desired length of the first vector $\vec{v_A}$ and finds the appropriate $\alpha$ value to meet it.

\section{Conclusion} \label{sec:conclusion}

We have presented an algorithm to interactively draw an extended log-aesthetic curve segment with a minor boundary condition by specifying the endpoints, the tangent vector at the first point, and a tangent line at the last point. The algorithm determines also the shape parameter of the log-aesthetic curve that eases the design since the control of the curve is based only on geometric data and no further parameter decision is required.

A log-aesthetic curve segment can be computed practically in real time, within milliseconds with a maximum error of $2\times10^{-10}$ \cite{yoshida2006interactive}. Our upgraded $\Lambda$ bisection algorithm extends the previous method only with conditions that do not increase the computation time significantly. Although, it is embedded into another $\alpha$ bisection that multiplies the running time, our experimental results, using an Intel i7 7700HQ CPU, shows that it still gives acceptable interactive real-time control. In our implementation, we calculated the curve segments using numerical integration. However, the precision can be increased and the computation time can be decreased (up to 13 times) using incomplete gamma functions \cite{ziatdinov2012analytic}.

\begin{figure}[htb]
   \centering
   \includegraphics[height=0.5\textheight]{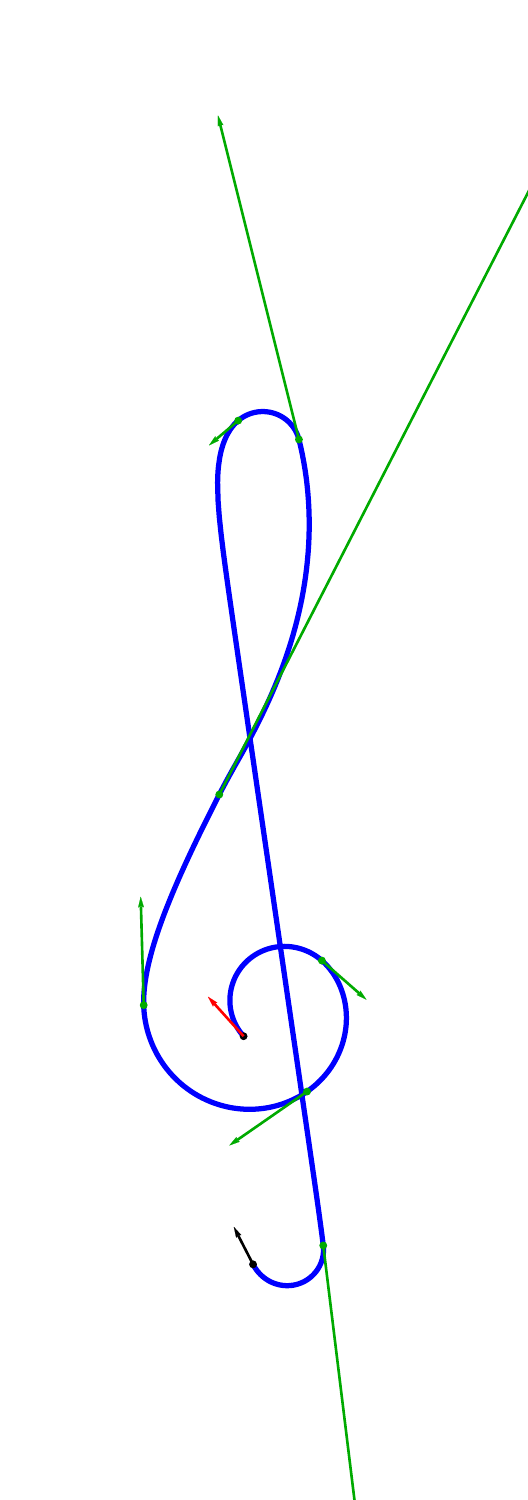}
    \caption{The violin (G) clef designed using the presented algorithm. The extended log-aesthetic curve segments are connected at the green points with $G^2$ continuity (both tangents and curvatures are match). The curve segments are controlled by the position of the points, by the direction of the tangent vectors (black, green, and red arrows) and by setting the length of the first tangent vector (red).
    }
\label{fig:violin}
\end{figure}

The method can also be used to $G^2$ connection of log-aesthetic curves effectively, that makes it more applicable for various tasks in the field of aesthetic design. An example is seen in \figref{fig:violin}. Regarding \cite{levien2009spiral}, to design fonts using aesthetic curve yields better results than using standard free-form curves because the design of font variation is more accessible and productive with these curves (e.g. interpolation between cubic B\'ezier-curves may fail to preserve even $G^1$-continuity). Our method could provide more intuitive design process with these curves.

\section*{Acknowledgements}
The first author was supported by the construction EFOP-3.6.3-VEKOP-16-2017-00002. The project was supported by the European Union, co-financed by the European Social Fund.

\appendix
\setcounter{secnumdepth}{0}
\section*{Pseudo-codes}

\begin{lstlisting}[
    caption={The new $\Lambda$ bisection algorithm (the angles are computed only as dot products). The modifications regarding the previous algorithm of \cite{yoshida2006interactive} are highlighted by red.},
    label={lst:lambda_bisection}
]
double lambda_bisection( double alpha, int max_iteration ) {
    double lmin = 0, lmax = 1, lambda, diff, theta_compare;
    int i = 0;
    bool enlarge, @beyond_inf_point = false@;
    @double point[2];  // 2D point to represent A_dash or B_dash @
    @if ( ( alpha <= 1 && !swap_flag ) || ( alpha >= 1 && swap_flag ) ) {
        theta_compare = calculate_theta_A_dash(lambda);
    } else {
        theta_compare = calculate_theta_C_dash(lambda);
    }@
    if ( alpha == 1 ) enlarge = true;
    else if ( alpha <= 1 ) {
        lmax = 1 / ( theta_$\Delta$ * (1-alpha) );
        @if ( theta_A < calculate_theta_A_dash(lmax) ) beyond_inf_point = true;@
        @point[2] = calculate_intersection_point(theta$\color{red}_\Delta$); // B_dash@
    } else {
        @lmax = 2 / ( -theta$\color{red}_\Delta$ * (1-alpha) );@
        @point[2] = calculate_curve_point(-theta$\color{red}_\Delta$); // A_dash@
    }
    lambda = ( lmin + lmax ) * 0.5;
    do {
        if ( alpha <= 1 ) diff = theta_compare - calculate_theta_A_dash(lambda);
        else diff = theta_compare - calculate_theta_C_dash(lambda);
        if ( diff < EPS ) return lambda; // success
        if ( ( 0 <= alpha && alpha <= 1 ) || 
             ( alpha < 1 && !beyond_inf_point ) ) {
            if ( diff < 0 ) {
                if ( enlarge ) lmax = lmax * 10;
                lmin = lambda;
                lambda += (lmax - lambda) * 0.5;
            } else {
                enlarge = false;
                lmax = lambda;
                lambda -= (lambda - lmin) * 0.5;
            }
        } @else if ( alpha < 1 && beyond_inf_point ) {
            if ( diff < 0 && point[0] > 0 ) {
                lmax = lambda;
                lambda -= (lambda - lmin) * 0.5;
            } else {
                lmin = lambda
                lambda += (lmax - lambda) * 0.5;
            }
        } else { // ( 1 < alpha )
            if ( diff < 0 || point[1] < 0 ) {
                lmax = lambda;
                lambda -= (lambda - lmin) * 0.5;
            } else {
                lmin = lambda;
                lambda += (lmax - lambda) * 0.5;
            }
        }@
        ++i;
    } while ( i < max_iteration );
    return -1; // not found
}
\end{lstlisting}

\begin{lstlisting}[
    caption={The $\alpha$ bisection algorithm that includes the $\Lambda$ bisection.},
    label={lst:alpha_bisection}
]
double alpha_bisection( double length, int max_iteration ) {
    double amin = -999, amax = 999; // arbitrary large interval
    double alpha = ( amin + amax ) * 0.5;
    int i = 0;
    do {
        lambda_bisection( alpha, max_iteration );
        double diff = length - Calculate_actual_length();
        if ( diff < EPS ) return alpha; // success
        if ( !swap_flag ) { // coordinates of the endpoints are not swapped
            if ( diff < 0 ) {
                amin = alpha;
                alpha += ( ( amax - alpha ) * 0.5 );
            } else {
                amax = alpha;
                alpha -= ( ( alpha - amin ) * 0.5 );
            }
        } else { // coordinates of the endpoints are swapped
            if ( instance_1 ) { // to include inflection point
                if ( alpha < 1 + EPSILON && beyond_inf_point ) {
                    // reversed sub-interval selection
                    if ( diff > 0 ) { 
                        amin = alpha;
                        alpha += ( ( amax - alpha ) * 0.5 );
                    } else {
                        amax = alpha;
                        alpha -= ( ( alpha - amin ) * 0.5 );
                    } continue;
                } else { // skipping unnecessary sub-interval
                    amax = alpha;
                    alpha -= ( ( alpha - amin ) * 0.5 );
                    continue;
                }
            }
            if ( instance_2 ) { // to include cusp point
                if ( alpha < 1 + EPSILON && beyond_inf_point ) {
                    // skipping unnecessary sub-interval
                    amin = alpha;
                    alpha += ( ( amax - alpha ) * 0.5 );
                    continue;
                }
            }
        }
        ++i;
    } while ( i < max_iteration );
    return 0; // not found
}
\end{lstlisting}

\bibliographystyle{plain}
\bibliography{main}

\end{document}